\documentstyle[12pt,fleqn,epsf]{article}
\newtheorem{theorem}{Theorem}
\newtheorem{proposition}{Proposition}
\newtheorem{lemma}{Lemma}
\newtheorem{definition}{Definition}
\newtheorem{conjecture}{Conjecture}

\renewcommand{\arraystretch}{1.3}
  
\def\nsection#1{\setcounter{equation}{0}\section{#1}}

\def\Bins#1#2{\left[{#1 \textstyle \atop #2}\right]}
\def\Bin#1#2{\biggl[{#1 \atop #2}\biggr]}
\def\Mults#1#2#3#4{\left[{#1 \textstyle \atop #2}\right]^{(#3)}_{#4}}
\def\Mult#1#2#3#4{\biggl[{#1 \atop #2}\biggr]^{(#3)}_{#4}}
\def\dMults#1#2#3#4{\left\{{#1 \textstyle \atop #2}\right\}^{(#3)}_{#4}}
\def\dMult#1#2#3#4{\biggl\{{#1 \atop #2}\biggr\}^{(#3)}_{#4}}
\def\mults#1#2{\left({\textstyle {#1 \atop #2} } \right)}
\def\mult#1#2{\biggl({#1 \atop #2}\biggr)}
\def\e{\mbox{e}}
\def\lf{\lfloor}
\def\rf{\rfloor}
\def\es{\mbox{\scriptsize e}}
\def\case#1#2{{\textstyle{#1\over #2}}}
\def\eps{\varepsilon}

\def\mod#1{\; (\bmod \: #1)}
\def\sc{\scriptstyle}

\begin{document}

\title{The Andrews-Gordon identities \\
and \\
$q$-multinomial coefficients}
 
\author{S.~Ole Warnaar\thanks{
e-mail: {\tt warnaar@maths.mu.oz.au}}
\\
Mathematics Department\\
University of Melbourne\\
Parkville, Victoria 3052\\
Australia}
  
\date{January, 1996 \\ \hspace{1mm}
 \\
Preprint No. 01-96}
   
\maketitle
   
\begin{abstract}
We prove polynomial boson-fermion
identities for the generating function
of the number of partitions of $n$ of the form 
$n=\sum_{j=1}^{L-1} j f_j$,
with $f_1\leq i-1$, $f_{L-1} \leq i'-1$ and
$f_j+f_{j+1}\leq k$.
The bosonic side of the identities involves $q$-deformations
of the coefficients of $x^a$ in the expansion of
$(1+x+\cdots+ x^k)^L$.
A combinatorial interpretation for these $q$-multinomial
coefficients is given using Durfee dissection partitions.
The fermionic side of the polynomial identities 
arises as the partition function of a one-dimensional
lattice-gas of fermionic particles.

In the limit $L\to\infty$, our identities reproduce the
analytic form of Gordon's generalization of the Rogers--Ramanujan
identities, as found by Andrews.
Using the $q \to 1/q$ duality, identities are obtained
for branching functions corresponding to cosets of type
$({\rm A}^{(1)}_1)_k \times ({\rm A}^{(1)}_1)_{\ell} /
({\rm A}^{(1)}_1)_{k+\ell}$ of fractional level $\ell$.
\end{abstract}
     
\newpage
      
\nsection{Introduction}\label{intro}
The Rogers--Ramanujan identities can be stated
as the following $q$-series identities.
\begin{theorem}[Rogers--Ramanujan]
For $a=0,1$ and $|q|<1$,
\begin{equation}
\sum_{n = 0}^{\infty} \frac{q^{n(n+a)}}{(1-q)(1-q^2)\ldots (1-q^n)}
=  \prod_{j=0}^{\infty}
(1-q^{5j+1+a})^{-1}
(1-q^{5j+4-a})^{-1}.
\label{RR}
\end{equation}
\end{theorem}
Since their independent discovery
by Rogers~[1-3], Ramanujan~\cite{Ramanujan} and 
also Schur~\cite{Schur}, many beautiful generalizations have been
found, mostly arising from partition-theoretic or Lie-algebraic
considerations, see refs.~\cite{Andrews76,LP} and references
therein. 

Most surprising, in 1981 Baxter rediscovered the Rogers--Ramanujan
identities (\ref{RR}) in his calculation 
of the order parameters of the hard-hexagon
model~\cite{Baxter81}, a lattice gas of hard-core particles
of interest in statistical mechanics.
It took however another ten years to fully realize the power 
of the (solvable) lattice model approach to finding $q$-series identities.
In particular, based on a numerical study of the
eigenspectrum of the critical three-state Potts model~\cite{KM,DKMM}
(yet another lattice model in statistical mechanics),
the Stony Brook group found an amazing variety of new $q$-series identities 
of Rogers--Ramanujan type~\cite{KKMMa,KKMMb}.
Almost none of these identities had been encountered
previously in the context of either partition theory
or the theory of infinite dimensional Lie algebras.

More specific, in the work of refs.~\cite{KKMMa,KKMMb}
expressions for Virasoro characters were given through
systems of fermionic quasi-particles. Equating these
{\em fermionic} character forms with the well-known 
Rocha-Caridi type {\em bosonic} expressions~\cite{RochaCaridi},
led to many $q$-series identities for Virasoro characters,
generalizing  the Rogers--Ramanujan identities 
(which are associated to the $M(2,5)$ minimal model).

\vspace*{5mm}
The proof of the Rogers--Ramanujan identities by means
of an extension to polynomial identities whose degree is determined by
a fixed integer $L$, was initiated by Schur~\cite{Schur}.
Before we elaborate on this approach, we need the 
combinatorial version of the Rogers--Ramanujan identities
stating that
\begin{theorem}[Rogers--Ramanujan]\label{RRC}
For $a=0,1$,
the partitions of $n$ into parts congruent to $1+a$ or $4-a$
$\mod{5}$ are equinumerous with the partitions of $n$ in
which the difference between any two parts is at least 2 and
1 occurs at most $1-a$ times.
\end{theorem}
Denoting the number of occurences of the part $j$ in a partition
by $f_j$, the second type of partitions in the above theorem
are those partitions
of $n=\sum_{\j \geq 1} j f_j$
which satisfy the following {\em frequency conditions}:
\begin{equation}
f_j+f_{j+1} \leq 1 \; \; \forall j
\qquad {\rm and} \qquad f_1 \leq 1-a.
\end{equation}
Schur notes that imposing the additional condition
$f_j=0$ for $j \geq L+1$, the generating
function of the ``frequency partitions'' satisfies the
recurrence 
\begin{equation}
g_L =  g_{L-1} + q^L g_{L-2}.
\label{Schurr}
\end{equation}
Together with the appropriate initial conditions,
Schur was able to solve these recurrences,
to obtain an alternating-sign type solution, now called a
bosonic expression.
Taking $L\to\infty$ in these bosonic polynomials yields
(after use of Jacobi's triple product identity) the
right-hand side of (\ref{RR}). Since this indeed
corresponds to the generating function of the ``$\mod{5}$''
partitions, this proves theorem~\ref{RRC}.
Much later, Andrews~\cite{Andrews70} obtained a solution to the
recurrence relation as a finite 
$q$-series with manifestly positive integer coefficients, 
now called a fermionic expression.
Taking $L\to\infty$ in these fermionic polynomials yields
the left-hand side of (\ref{RR}).

Recently much progress has been made in proving
the boson-fermion identities of~\cite{KKMMa,KKMMb}
(and generalizations thereof), by following the Andrews--Schur
approach.
That is, for many of the Virasoro-character identities, finitizations 
to polynomial boson-fermion identities have been found, which
could then be proven either fully recursively (\`a la Andrews)
or one side combinatorially and one side 
recursively (\`a la Schur), see refs.~[15--30].

\vspace*{5mm}
In this paper we consider polynomial identities which imply
the Andrews--Gordon generalization of the Rogers--Ramanujan identities.
First, Gordon's theorem~\cite{Gordon},
which provides a combinatorial generalization
of the Rogers--Ramanujan identities, reads
\begin{theorem}[Gordon]\label{GI}
For all $k\geq 1$, $1\leq i\leq k+1$,
let $A_{k,i}(n)$ be the number of partitions of $n$
into parts not congruent to $0$ or $\pm i \mod{2k+3}$
and let $B_{k,i}(n)$ be the number of
partitions of $n$ of the form
$n = \sum_{j\geq 1} j f_j$, with $f_1\leq i-1$ and
$f_j+f_{j+1}\leq k$ (for all $j$). 
Then $A_{k,i}(n)=B_{k,i}(n)$.
\end{theorem}
Subsequently the following analytic counterpart of this result 
was obtained by Andrews~\cite{Andrews74}, 
generalizing the analytic form (\ref{RR}) of the
Rogers--Ramanujan identities.
\begin{theorem}[Andrews]
For all $k\geq 1$, $1\leq i\leq k+1$ and $|q|<1$,
\begin{equation}
\renewcommand{\arraystretch}{0.8}
\sum_{n_1,n_2,\ldots,n_{k}\geq 0}
\frac{ q^{N_1^2 + \cdots + N_k^2 + N_i + \cdots + N_k}}
{(q)_{n_1} (q)_{n_2} \cdots (q)_{n_k} } 
=  \prod_{
\begin{array}{c}
\sc j=1 \\
\sc j \not\equiv 0,\pm i \mod{2k+3}
\end{array}}^{\infty} (1-q^j)^{-1}
\label{An}
\end{equation}
with
\begin{equation}
N_j = n_j + \cdots + n_{k}
\label{Nj}
\end{equation}
and $(q)_a=\prod_{k=1}^a (1-q^k)$ for $a>0$ and $(q)_0=1$.
\end{theorem}
Application of Jacobi's triple product identity admits for a rewriting
of the right-hand side of (\ref{An}) to
\begin{equation}
\frac{1}{(q)_{\infty}}
\sum_{j=-\infty}^{\infty} (-)^j 
q^{j\bigl( (2k+3)(j+1)-2i\bigr)/2}.
\label{Anb}
\end{equation}

Equating (\ref{Anb}) and the left-hand side of (\ref{An}),
gives an example of a boson-fermion identity.
Here we consider, in the spirit of Schur,
a ``natural'' finitization of Gordon's 
frequency condition such that this boson-fermion identity 
is a limiting case of polynomial identities.
In particular, we are interested in the quantity
$B_{k,i,i';L}(n)$, counting the number of partitions
of $n$ of the form
\begin{equation}
n = \sum_{j=1}^{L-1} j f_j
\end{equation}
with frequency conditions
\begin{equation}
f_1\leq i-1, \qquad f_{L-1} \leq i'-1 \qquad {\rm and}
\qquad 
f_j+f_{j+1}\leq k \quad \mbox{for $j=1,\ldots,L-2$}.
\end{equation}
If we denote the generating function of
partitions counted by $B_{k,i,i';L}(n)$ by
$G_{k,i,i';L}(q)$, then clearly
$\lim_{L\to\infty}
G_{k,i,i';L}(q)=G_{k,i}(q)$, with
$G_{k,i}$ the generating function
associated with $B_{k,i}(n)$ of theorem~\ref{GI}.
Also note that $G_{k,i,1;L}=G_{k,i,k+1;L-1}$.
Our main results can be formulated as the following
two theorems for $G_{k,i,i';L}$.

Let $\Bin{L}{a}$ be the Gaussian polynomial or $q$-binomial coefficient
defined by
\begin{equation}
\Bin{L}{a}= \Bin{L}{a}_q = \left\{
\begin{array}{ll}
\displaystyle \frac{(q)_L}
{(q)_{a}(q)_{L-a}} \qquad & 0\leq a \leq L \\[3mm]
0& \mbox{otherwise.}
\end{array}
\right.
\label{Gpoly}
\end{equation}
Further, let ${\cal I}_k$ be the incidence matrix
of the Dynkin diagram of A$_k$ with an additional
tadpole at the $k$-th node: 
\begin{equation}
({\cal I}_k)_{j,\ell}
=\delta_{j,\ell-1}+\delta_{j,\ell+1}+
\delta_{j,\ell} \delta_{j,k} \qquad
\qquad j,\ell=1,\ldots,k,
\end{equation}
and let $C_k$ be the corresponding Cartan-type
matrix, $(C_k)_{j,\ell}=2 \delta_{j,\ell}-
({\cal I}_k)_{j,\ell}$.
Finally let $\vec{n}$, $\vec{m}$ and $\vec{\e}_j$ be
$k$-dimensional (column)-vectors with entries
$\vec{n}_j=n_j$, $\vec{m}_j=m_j$
and $(\vec{\e}_j)_{\ell}=
\delta_{j,\ell}$. Then
\begin{theorem}\label{ft}
For all $k\geq 1$, $1\leq i,i'\leq k+1$ and
$kL\geq 2k-i-i'+2$, 
{\rm 
\begin{equation}
G_{k,i,i';L}(q)=
\sum_{n_1,n_2,\ldots,n_k\geq 0}
q^{\displaystyle \, \vec{n}^T C^{-1}_k (\vec{n} 
+\vec{\e}_k -\vec{\e}_{i-1})}
\prod_{j=1}^k
\Bin{n_j+m_j}{n_j},
\label{Fermion}
\end{equation}}
with $(m,n)$-system \cite{Berkovich} given by
{\rm 
\begin{equation}
\vec{m}+\vec{n}=
\case{1}{2}\Bigl({\cal I}_k \, \vec{m}
+(L\!-\! 2) \, \vec{\e}_k + \vec{\e}_{i-1}
+\vec{\e}_{i'-1}\Bigr).
\label{mn}
\end{equation}}
\end{theorem}
We note that $(C^{-1}_k)_{j,\ell}=\min(j,\ell)$
and hence that, using the variables $N_j$ of 
(\ref{Nj}), we can rewrite the quadratic
exponent of $q$ in (\ref{Fermion}) as
$N_1^2+ \cdots + N_k^2+N_i+\cdots+ N_k$.
For $k\geq 2$, the ``finitization'' (\ref{Fermion})-(\ref{mn}) of the
left-hand side of (\ref{An}) is new. 
For $k=1$ it is the already mentioned
fermionic solution to the
recurrence (\ref{Schurr}) as found by Andrews~\cite{Andrews70}. 
Another finitization,
which does not seem to be related to a finitization
of Gordon's frequency conditions,
has recently been proposed in refs.~[17-19]
(see also ref.~\cite{BM94}).
A more general expression, which
includes (\ref{Fermion})-(\ref{mn}) and that of [17-19] as special cases,
will be discused in section~\ref{secdis}.

Our second result, which is maybe of more interest
mathematically since it involves new generalizations of the
Gaussian polynomials, can be stated as follows.
Let $\Mults{L}{a}{p}{k}$ be the $q$-multinomial
coefficient defined in equation (\ref{qM}) of the 
subsequent section. 
Also, define 
\begin{equation}
r=k-i'+1
\label{defr}
\end{equation}
and 
\begin{equation}
s = \left\{\begin{array}{lcl}
i \qquad & \mbox{for} & i=1,3,\ldots,2\lfloor\frac{k}{2}\rfloor +1 
\\[1mm]
2k+3-i & \mbox{for} & i=2,4,\ldots,2 \lfloor\frac{k+1}{2} \rfloor,
\end{array} \right.
\label{defs}
\end{equation}
so that $r=0,1,\ldots,k$ and $s=1,3,\ldots,2k+1$.
Then
\begin{theorem}\label{bt}
For all $k>0$, $1\leq i,i'\leq k+1$ and
$kL\geq 2k-i-i'+2$, 
\begin{eqnarray}
G_{k,i,i';L}(q)&=&
\sum_{j=-\infty}^{\infty} \left\{
q^{j\bigl((2j+1)(2k+3)-2s\bigr)}
\Mult{L}{\case{1}{2}(k L+k-s-r+1)+(2k+3)j}{r}{k} 
\right.  \nonumber \\
& & \hspace{10mm} \left.
-q^{\bigl(2j+1\bigr)\bigl((2k+3)j+s\bigr)}
\Mult{L}{\case{1}{2}(k L+k+s-r+1)+(2k+3)j}{r}{k}
\right\}
\label{Boson1}
\end{eqnarray}
for $r\equiv k(L+1) \mod{2}$ and
\begin{eqnarray}
G_{k,i,i';L}(q)&=&
\sum_{j=-\infty}^{\infty} \left\{
q^{j\bigl((2j+1)(2k+3)-2s\bigr)}
\Mult{L}{\case{1}{2}(k L-k+s-r-2)-(2k+3)j}{r}{k} 
\right.  \nonumber \\
& & \hspace{10mm} \left.
-q^{\bigl(2j+1\bigr)\bigl((2k+3)j+s\bigr)}
\Mult{L}{\case{1}{2}(k L-k-s-r-2)-(2k+3)j}{r}{k}
\right\}
\label{Boson2}
\end{eqnarray}
for $r\not\equiv k(L+1) \mod{2}$.
\end{theorem}
For $k\geq 3$,
the finitizations (\ref{Boson1}) and (\ref{Boson2}) of the
right-hand side of (\ref{An}) are new.
For $k=1$ (\ref{Boson1}) and (\ref{Boson2}) are Schur's bosonic polynomials.
For $k=2$, $\Mults{L}{a}{p}{2}$ being a $q$-trinomial coefficient, 
(\ref{Boson1}) and (\ref{Boson2}) were (in a slightly different
representation) first obtained in ref.~\cite{AB}.
An altogether different alternating-sign expression for $G_{k,i,i';L}$
in terms of $q$-binomials has been found in ref.~\cite{DJKMO}.
A different finitization of the right-hand side
of (\ref{An}) involving $q$-binomials has been given 
in~\cite{Andrews70,ABF}.
A more general expression, which includes (\ref{Boson1}), (\ref{Boson2}) 
and that of ref.~\cite{Andrews70,ABF} as special cases,
will be discused in section~\ref{secdis}.

Equating (\ref{Fermion}) and (\ref{Boson1})--(\ref{Boson2})
leads to non-trivial polynomial identities, which in the
limit $L\to\infty$ reduce to Andrews' analytic form of
Gordon's identity. For $k=1$ these are the polynomial identities
featuring in the Andrews--Schur proof
of the Rogers--Ramanujan identities (\ref{RR})~\cite{Andrews70}.

\vspace*{5mm}
The remainder of the paper is organized
as follows. In the next section we
introduce the $q$-multinomial coefficients
and list some $q$-multinomial identities needed
for the proof of theorem~\ref{bt}.
Then, in section~\ref{secC}, a combinatorial interpretation
of the $q$-multinomials is given using Andrews' Durfee dissection
partitions. In section~\ref{pt4} we give a recursive 
proof of theorem~\ref{bt}
and in section~\ref{pt3} we prove theorem~\ref{ft} combinatorially, 
interpreting the restricted frequency partitions as configurations of a 
one-dimensional lattice-gas of fermionic particles.
We conclude this paper with a discussion of our results,
a conjecture generalizing theorems~\ref{ft} and~\ref{bt},
and some new identities for the branching functions of
cosets of type
$({\rm A}^{(1)}_1)_k \times ({\rm A}^{(1)}_1)_{\ell} /
({\rm A}^{(1)}_1)_{k+\ell}$ with fractional level $\ell$.
Finally, proofs of some of the $q$-multinomial identities 
are given in the appendix.

\nsection{$q$-multinomial coefficients}\label{secqB}
Before introducing the $q$-multinomial coefficients,
we first recall some
facts about ordinary multinomials.
Following ref.~\cite{AB}, we define $\mults{L}{a}_k$ 
for $a=0,\ldots,kL$ as
\begin{equation}
(1+x+\cdots + x^k)^L = \sum_{a=0}^{k L} \:
\mult{L}{a}_k \: x^a.
\end{equation}
Multiple use of the binomial theorem yields
\begin{equation}
\mult{L}{a}_k =
\sum_{j_1+ \cdots + j_k = a}
\mult{L}{j_1}
\mult{j_1}{j_2}
\cdots
\mult{j_{k-1}}{j_k},
\label{multq1}
\end{equation}
where $\mults{L}{a}=\mults{L}{a}_1$ is the usual
binomial coefficient.
 
Some readily established properties of $\mults{L}{a}_k$ are
the symmetry relation
\begin{equation}
\mult{L}{a}_k = \mult{L}{kL-a}_k
\label{sym}
\end{equation}
and the recurrence
\begin{equation}
\mult{L}{a}_k =
\sum_{m=0}^k
\mult{L-1}{a-m}_k .
\label{rec}
\end{equation}
  
For our subsequent working it will be convenient to define $k+1$
different $q$-deformations of the
multinomial coefficient (\ref{multq1}).
\begin{definition}
For $p=0,\ldots,k$ we set
\begin{equation}
\Mult{L}{a}{p}{k} =
\sum_{j_1+ \cdots + j_k = a}
q^{\;\displaystyle
\sum_{\ell=1}^{k-1} (L-j_{\ell})j_{\ell+1}-
\sum_{\ell=k-p}^{k-1} j_{\ell+1}}
\Bin{L}{j_1}
\Bin{j_1}{j_2}
\cdots
\Bin{j_{k-1}}{j_k},
\label{qM}
\end{equation}
with $\Bins{L}{a}$ the standard $q$-binomial coefficients
of (\ref{Gpoly}).
\end{definition}
Note that $\Mults{L}{a}{p}{k}$ is unequal to zero for 
$a=0,\ldots,kL$ only. Also note the initial condition
\begin{equation}
\Mults{0}{a}{p}{k}  = \delta_{a,0}.
\label{Lis0}
\end{equation}

In the following we state a number of $q$-deformations
to (\ref{sym}) and (\ref{rec}).
Although our list is certainly not exhaustive, we
have restricted ourselves to those identities
which in our view are simplest, and to those 
needed for proving theorem~\ref{bt}.
Most of these identities are generalizations of known
$q$-binomial and $q$-trinomial identities which, for example,
can be found in refs.~\cite{AB,Andrews94,BMO,BM95}.

First we put some simple symmetry properties
generalizing (\ref{sym}), in a lemma.
\begin{lemma}\label{lemsym}
For $p=0,\ldots,k$ the following symmetries hold:
\begin{equation}
\Mult{L}{a}{p}{k} = q^{(k-p)L-a} 
\Mult{L}{kL-a}{k-p}{k} 
\qquad {\mbox and}
\qquad
\Mult{L}{a}{0}{k} = 
\Mult{L}{kL-a}{0}{k}.
\label{qs}
\end{equation}
\end{lemma}
The proof of this lemma is given in the appendix.

To our mind the simplest way of $q$-deforming
(\ref{rec}) (which was communicated to us by
A. Schilling) is
\begin{proposition}[Fundamental recurrences; Schilling]\label{pfr}
For $p=0,\ldots,k$, the $q$-multinomials
satisfy
\begin{equation}
\Mult{L}{a}{p}{k}=
\sum_{m=0}^{k-p} q^{m(L-1)} \Mult{L-1}{a-m}{m}{k}+
\sum_{m=k-p+1}^k q^{L(k-p)-m} \Mult{L-1}{a-m}{m}{k}.
\label{frec}
\end{equation}
\end{proposition}
In the next section we give a combinatorial proof of this important
result for the $p=0$ case.
An analytic proof for general $p$ has been given by Schilling
in ref.~\cite{Schillingb}.

We now give some equations, proven in the appendix,
which all reduce to the tautology $1=1$ in the $q\to 1$ limit.
\begin{proposition}\label{ptaut}
For all $p=-1,\ldots,k-1$, we have
\begin{equation}
\Mult{L}{a}{p}{k}
+ q^L \Mult{L}{kL-a-p-1}{p+1}{k}
=
\Mult{L}{kL-a-p-1}{p}{k}
+q^L \Mult{L}{a}{p+1}{k} ,
\label{taut}
\end{equation}
with $\Mults{L}{a}{-1}{k}=0$.
\end{proposition}
The power of these ($q$-deformed) tautologies is that they allow 
for an endless number of different rewritings of the fundamental 
recurrences.
In particular, as shown in the appendix, they allow for the 
non-trivial transformation of (\ref{frec}) into
\begin{proposition}\label{pqr2}
For all $p=0,\ldots,k$, we have
\begin{eqnarray}
\Mult{L}{a}{p}{k} &=&
\renewcommand{\arraystretch}{0.6}
\sum_{\begin{array}{c}
\sc m=0 \\
\sc m\equiv p+k \mod{2}
\end{array} }^{k-p}
q^{m(L-1)}
\Mult{L-1}{a-\case{1}{2}(m-p+k)}{m}{k}
\nonumber \\
&+&
\renewcommand{\arraystretch}{0.6}
\sum_{\begin{array}{c}
\sc m=0 \\
\sc m\not\equiv p+k \mod{2}
\end{array} }^{k-p-1}
q^{m(L-1)}
\Mult{L-1}{kL-a-\case{1}{2}(m+p+k+1)}{m}{k}
\nonumber \\
&+&
\renewcommand{\arraystretch}{0.6}
\sum_{\begin{array}{c}
\sc m=k-p+2 \\
\sc m\equiv p+k \mod{2}
\end{array} }^k
q^{\case{1}{2}\bigl((2L-1)(k-p)-m\bigr)}
\Mult{L-1}{a-\case{1}{2}(m-p+k)}{m}{k}
\nonumber \\
&+&
\renewcommand{\arraystretch}{0.6}
\sum_{\begin{array}{c}
\sc m=k-p+1 \\
\sc m\not\equiv p+k \mod{2}
\end{array} }^k
q^{kL+\case{1}{2}\bigl((2L+1)(k-p)-m+1\bigr)-2a}
\Mult{L-1}{kL-a-\case{1}{2}(m+p-k-1)}{m}{k}.
\label{qrec}
\end{eqnarray}
\end{proposition}
It is thanks to these rather unappealing
recurrences that we can prove theorem~\ref{bt}.

Before concluding this section on the $q$-multinomial
coefficients let us make some further remarks.
First, for $k=1$ and $k=2$ we reproduce the well-known
$q$-binomial and $q$-trinomial coefficients.
In particular, 
\begin{equation}
\Mult{L}{a}{0}{1} = \Bin{L}{a}
\end{equation}
and 
\begin{equation}
\Mult{L}{a}{p}{2} = \mult{L;L-a-p;q}{L-a}_2 \qquad
\mbox{ for } p=0,1,
\label{trin}
\end{equation}
where on the right-hand side of (\ref{trin}) we have
used the $q$-trinomial notation introduced
by Andrews and Baxter~\cite{AB}.

Second, in \cite{AB}, several
recurrences involving $q$-trinomials with
just a single superscript $(p)$ are given.
We note that such recurrences follow from (\ref{frec})
by taking the difference between various values
of $p$.
In particular we have for all $r=0,\ldots,p$
\begin{eqnarray}
\Mult{L}{a}{p}{k} = \Mult{L}{a}{p-r}{k}
&+& q^{L(k-p)-a} \sum_{m=0}^{p-r-1} \left(1-q^{rL}\right)
q^{m(L-1)} \Mult{L-1}{kL-a-m}{m}{k}
\nonumber \\
&+& q^{L(k-p)-a} \sum_{m=p-r}^{p-1} \left(1-q^{(p-m)L}\right)
q^{m(L-1)} \Mult{L-1}{kL-a-m}{m}{k}.
\label{pr}
\end{eqnarray}
This can be used to eliminate all multinomials
$\Mults{..}{..}{m}{k}$ for $m=0,\ldots,p-1,p+1,\ldots,k$
in favour of $\Mults{..}{..}{p}{k}$. 
The price to be paid for this is that the resulting
expressions tend to get very complicated if $k$ gets large.

A further remark we wish to make is that to 
our knowledge the general $q$-deformed multinomials
as presented in (\ref{qM}) are new.
The multinomial $\Mults{L}{a}{0}{k}$ however was already
suggested as a ``good'' $q$-multinomial by Andrews in
\cite{Andrews94}, where the following 
generating function for $q$-multinomials was proposed
for all $k>1$:
\begin{equation}
p_{k,L}(x) = \sum_{a=0}^L x^a q^{\mults{a}{2}} 
\Bin{L}{a} \; p_{k-1,a}(xq^L),
\end{equation}
with $p_{0,L}(x)=1$.
Clearly,
\begin{equation}
p_{k,L} = \sum_{a=0}^{kL} x^a q^{\mults{a}{2}} \Mult{L}{a}{0}{k}.
\end{equation}
Also in the work of Date {\em et al.} the $\Mult{L}{a}{0}{k}$
makes a brief appearance, see ref.~\cite{DJKMO} eqn.~(3.29).

The more general $q$-multinomials of equation (\ref{qM}) have been
introduced independently by Schilling~\cite{Schillingb}.
(The notation used in ref.~\cite{Schillingb} and that of the present
paper is almost identical apart from the fact that
$\Mults{L}{a}{p}{k}$ is replaced by
$\Mults{L}{kL/2-a}{p}{k}$.)

\nsection{Combinatorics of $q$-multinomial coefficients}\label{secC}
In this section a combinatorial
interpretation of the $q$-multinomials coefficients is given
using Andrews' Durfee dissections~\cite{Andrews79}. 
We then show how the fundamental recurrences
(\ref{frec}) with $p=0$ follow as an immediate consequence of this
interpretation.

As a first step it is convenient to change
variables from $q$ to $1/q$.
Using the elementary transformation property of the
Gaussian polynomials
\begin{equation}
\Bin{L}{a}_{1/q}  = 
q^{-a(L-a)}\Bin{L}{a}_{q}, 
\end{equation}
we set
\begin{definition}\label{dd}
For $p=0,\ldots,k$
{\rm 
\begin{eqnarray}
\dMult{L}{a}{p}{k} &:=&
\left. q^{-a L} \Mult{L}{a}{p}{k}  \; \right|_{q\to 1/q}
\nonumber \\
&=& 
\sum_{N_1+ \cdots + N_k = a}
q^{N_1^2 + \cdots + N_k^2
+N_{k-p+1}+ \cdots + N_k}
\Bin{L}{N_1}
\Bin{N_1}{N_2}
\cdots
\Bin{N_{k-1}}{N_k} 
\label{dqM} \\
&=& 
\sum_{\vec{n}^T C^{-1}_k \vec{\es}_k=a}
q^{\displaystyle \, \vec{n}^T C^{-1}_k (\vec{n} 
+\vec{\e}_k -\vec{\e}_{k-p})}
\frac{(q)_L}{(q)_{L-\vec{n}^T C^{-1}_k \vec{\es}_1} (q)_{n_1}
(q)_{n_2} \ldots (q)_{n_k}}\; .
\end{eqnarray}}
\end{definition}

\subsection{Successive Durfee squares and Durfee dissections}
As a short intermezzo, we review some of the ideas introduced by
Andrews in ref.~\cite{Andrews79}, needed for our interpretation of
(\ref{dqM}).
Those already familiar with such concepts as ``$(k,a)$-Durfee dissection
of a partition'' and ``$(k,a)$-admissible partitions'' may wish to skip
the following and resume in section~\ref{next}.
Throughout the following a partition and its corresponding
Ferrers graph are identified.

\begin{figure}[t]
\centerline{\epsffile{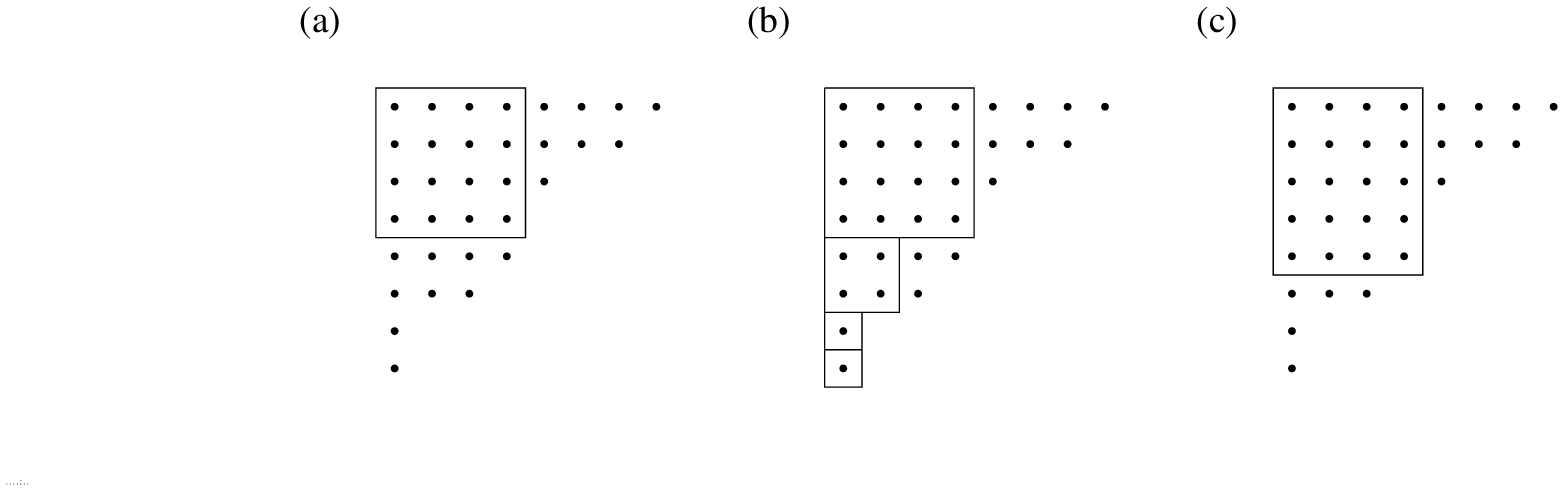}}
\caption{(a) Durfee square of the partition 
$\pi_{\rm ex}=9+7+5+4+4+3+1+1$. (b) The four successive Durfee
squares of $\pi_{\rm ex}$. (c) The Durfee rectangle of $\pi_{\rm ex}$.}
\label{DS}
\end{figure}

\begin{definition}
The Durfee square of a partition is the maximal square of nodes
(including the upper-leftmost node).
\end{definition}
The {\em size} of the Durfee square is the number of rows
for which $r_{\ell}\geq \ell$,
labelling the rows (=parts) of a partition by $r_1\geq r_2 \geq \ldots$.
Copying the example from ref.~\cite{Andrews79}, the Ferrers graph and
Durfee square of the partition $\pi_{\rm ex}=9+7+5+4+4+3+1+1$ is shown 
in Figure~\ref{DS}(a).

The portion of a partition of $n$ below its Durfee square defines 
a partition of $m<n$. For this ``smaller'' partition one can again draw
the Durfee square.
Continuing this process of drawing squares,
we end up with the successive Durfee squares of a partition.
For the partition $\pi_{\rm ex}$ this is shown in Figure~\ref{DS}(b).
If a partition $\pi$ has $k$ successive Durfee squares,
with $N_{\ell}$ the size of the $\ell$-th square,
then $\pi$ has exactly $N_1+\cdots + N_k$ parts with
$N_1+\cdots+N_{\ell}$ parts $\geq N_{\ell}$
for all $\ell=1,\ldots,k$. 

Following Andrews we now slightly generalize the previous notions.
\begin{definition}[Durfee rectangle]
The Durfee rectangle of a partition is the maximal rectangle
of nodes
whose height exceeds its width by precisely one row.
\end{definition}
The Durfee rectangle of the partition $\pi_{\rm ex}$ is shown in
Figure~\ref{DS}(c).
The {\em size} of the Durfee rectangle is its width.

One can now combine the Durfee squares and rectangles to define
\begin{definition}[Durfee dissection]
The $(k,i)$-Durfee dissection of a partition
is obtained by drawing $i-1$ successive Durfee
squares followed by $k-i+1$ successive Durfee rectangles.
\end{definition}
In the following it will be convenient to adopt a slightly
unconventional labelling. In particular,
we label the Durfee squares from 1 to $i-1$
and the rectangles from $i$ to $k$. Correspondingly,
$N_{\ell}$ is the size of the Durfee square or rectangle
labelled by $\ell$.
We note that in the $(k,i)$-dissection of a partition corresponding to
Durfee squares and rectangles of respective sizes
$N_1\geq N_2 \geq \ldots \geq N_k$, all  
the $N_{\ell}$ beyond some fixed $\ell'$ may actually be zero.

Finally we come to the most important definition of this section.
\begin{definition}[$(k,i)$-admissible]
Let $N_1\geq N_2 \geq \ldots \geq N_k$ be
the respective sizes of the Durfee squares and rectangles
in the $(k,i)$-Durfee dissection  of a partition $\pi$.
Then $\pi$  is $(k,i)$-admissible if
\begin{itemize}
\item
$\pi$ has no parts below its last successive Durfee
rectangle (or square if $i=k+1$.)
\item
For $\ell=i,\ldots,k$,
the last row of the Durfee rectangle labelled by $\ell$ has
$N_{\ell}$ nodes.
\end{itemize}
\end{definition}
The first condition is equivalent to stating
that the number of parts of $\pi$ equals
$N_1+\cdots + N_k + \max(\ell'-i+1,0)$,
where $\ell'$ labels the number of Durfee squares and rectangles
of non-zero size; $N_{\ell}>0$ for $\ell \leq \ell'$ and
$N_{\ell}=0$ for $\ell>\ell'$.
The second condition is equivalent to stating  that the last
row of each Durfee rectangle is actually a part of $\pi$.

\subsection{$(k,i;L,a)$-admissible partitions and $q$-multinomial
coefficients}\label{next}
Using the previous definitions we are now prepared for the
combinatorial interpretation of (\ref{dqM}).

\begin{figure}[t]
\epsfxsize = 13cm
\centerline{\epsffile{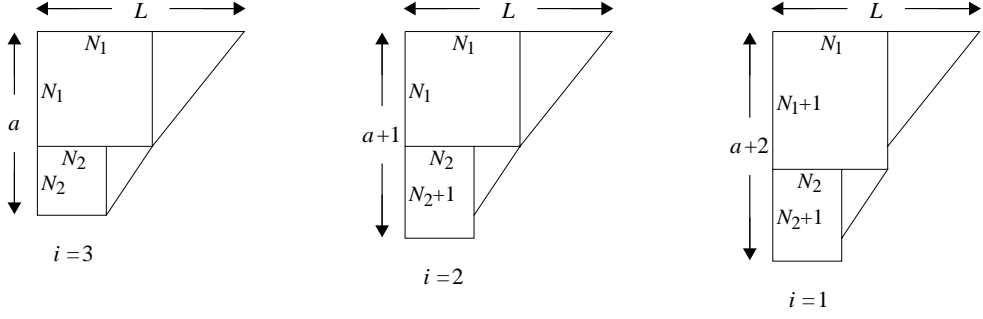}}
\caption{Graphical representation of 
the $(2,i;L,a)$-admissible partitions,
generated by $\dMults{L}{a}{3-i}{2}$.
The respective values of $N_1$ and $N_2$ are free to vary,
only their sum taken the fixed value $a$.
Note that the number of parts in the second and third figure
are actually not fixed, but vary between $a$ and $a-i+3$, depending
on the number of Durfee rectangles of non-zero size.}
\label{combin}
\end{figure}

\begin{definition}[$(k,i;L,a)$-admissible]
Let $N_1\geq N_2 \geq \ldots \geq N_k$ be
the respective sizes of the Durfee squares and rectangles
of a $(k,i)$-admissible partition $\pi$.
Then $\pi$ is said to be
$(k,i;L,a)$-admissible if 
the largest part of $\pi$ is less or equal to $L$ and
$N_1+\cdots + N_k=a$.
\end{definition}
For a $(k,i;L,a)$-admissible partition $\pi$, the portion $\pi_{\ell}$
of $\pi$ to the right of the Durfee square or rectangle labelled
by $\ell$ (and below the Durfee square or rectangle 
labelled $N_{\ell-1}$), 
is a partition with largest part $\leq N_{\ell-1}-N_{\ell}$
(where $N_0=L$) and number of parts $\leq N_{\ell}$.
Recalling that the Gaussian polynomial (\ref{Gpoly}) 
is the generating function of partitions with largest part $\leq L-a$ and
number of parts $\leq a$~\cite{Andrews76},
we thus find that the generating function of
$(k,i;L,a)$-admissible partitions is given by
\begin{equation}
\sum_{N_1+ \cdots + N_k = a}
q^{N_1^2} \Bin{L}{N_1} \cdots
q^{N_{i-1}^2} \Bin{N_{i-2}}{N_{i-1}} 
q^{N_i(N_i+1)} \Bin{N_{i-1}}{N_i} \cdots
q^{N_k(N_k+1)} \Bin{N_{k-1}}{N_k} = 
\dMult{L}{a}{k-i+1}{k}.
\end{equation}
Denoting an arbitrary partition 
of $n$ with largest part $\leq L$ and number of parts $\leq a$ by a
rectangle of width $L$ and height $a$, the $(k,i;L,a)$-admissible
partitions can be represented graphically as shown 
in Figure~\ref{combin} for the case $k=2$.

Equipped with the above interpretation we return to the
recurrence relation (\ref{frec}) for $p=0$.
Using definition~\ref{dd} to rewrite this in terms of 
$\dMults{L}{a}{p}{k}$, gives
\begin{equation}
\dMult{L}{a}{0}{k}=
q^a \sum_{m=0}^{k} \dMult{L-1}{a-m}{m}{k}.
\end{equation}
This is obviously true if the following combinatorial statements
hold.
\begin{lemma}\label{s2}
\hfill
\begin{itemize}
\item
Adding a column of $a$ nodes to the left 
of a $(k,k-m+1;L-1,a-m)$-admissible partition with
$m\in \{0,1,\ldots,k\}$, yields
a $(k,k+1;L,a)$-admissible partition.
\item
Removing the first column (of $a$ nodes) from a 
$(k,k+1;L,a)$-admissible partition yields a
$(k,k-m+1;L-1,a-m)$-admissible partition 
for some $m\in \{0,1,\ldots,k\}$.
\end{itemize}
\end{lemma}
To show the first statement,
we note that a partition is $(k,k+1;L,a)$-admissible
if it has exactly $a$ parts, has largest part $\leq L$ and has at most
$k$ successive Durfee squares.
A $(k,k-m+1;L-1,a-m)$-admissible partition  has at most $a$ parts and
has largest part $\leq L-1$. Hence adding a column of $a$ nodes
to the left of such a partition, yields a partition $\pi$
which has $a$ parts and largest part $\leq L$.
Remains to show that $\pi$ has at most
$k$ successive Durfee squares.
To see this first assume that 
the $(k,k-m+1;L-1,a-m)$-admissible partition only consists of
Durfee squares and rectangles. That is, we have a partition of
$N_1^2+\cdots + N_k^2 + N_{k-m+1} + \cdots + N_k$, with
$N_1+ \cdots + N_k=a-m$.
Adding a column of $a$ dots trivially yields a partition $\pi$
with $k$ successive Durfee squares with respective sizes
\begin{equation}
N_1 \geq N_2 \geq ... \geq N_{k-m} \geq N_{k-m+1}+1 \geq \ldots
\geq N_k+1>0,
\end{equation}
with $\pi$ having a column of $N_{\ell}$  
nodes to the right of the $\ell$-th successive Durfee square
for each $\ell \leq k-m$.
Now note that we in fact have treated the ``worst'' possible cases.
All other $(k,k-m+1;L-1,a-m)$-admissible partitions 
can be obtained from the ``bare'' ones just treated by adding
partitions with largest part 
$\leq N_{\ell-1}-N_{\ell}$
(where $N_0=L$) and number of parts $\leq N_{\ell}$
to the right of the Durfee square or rectangle labelled by $\ell$ for all
$\ell$.
Let $\pi$ be such a ``dressed'' partition, obtained from a bare
$(k,k-m+1;L-1,a-m)$-admissible partition $\pi_b$, and let the images
of $\pi$ and $\pi_b$ after adding a column of $a$ dots be $\pi'$ and $\pi_b'$.
Further, let $N_{\ell}$ and $M_{\ell}$ be the size of the 
$\ell$-th successive Durfee
square of $\pi_b'$ and $\pi'$, respectively. 
Since $\pi$ is obtained from $\pi_b$ by adding additional nodes
to its rows, we have $M_1+\cdots M_{\ell} \geq N_1+\cdots+N_{\ell}$
for all $\ell$. 
From the fact that $\pi_b'$ has at most $k$ successive Durfee squares
it thus follows that this is also true for $\pi'$.

To show the second statement of the lemma, note that 
from (\ref{rec}) we see that the map implied by
the first statement is in fact a map onto the set
of $(k,k+1,L,a)$-admissible partitions.
Since for $m\neq m'$,
the set of $(k,k-m+1;L-1,a-m)$-admissible partitions is distinct from
the set of $(k,k-m'+1;L-1,a-m')$-admissible partitions,
the second statement immediately follows.

To prove (\ref{frec}) is true for general $p$, we need to establish
\begin{equation}
\dMult{L}{a}{p}{k}=
q^a \sum_{m=0}^{k-p} \dMult{L-1}{a-m}{m}{k}
+q^a \sum_{m=k-p+1}^k q^{L(p-k+m)} \dMult{L-1}{a-m}{m}{k}.
\end{equation}
Unfortunately, a generalization of lemma~\ref{s2}
which would imply this more general result has so far eluded us.

Before concluding our discussion of $q$-multinomial
coefficients we note that if the restriction
on $L$ is dropped in the $(k,i;L,a)$-admissible
partitions, their generating function reduces to
\begin{equation}
\renewcommand{\arraystretch}{0.7}
\lim_{L\to\infty} \dMult{L}{a}{k-i+1}{k}=
\sum_{
\begin{array}{c}
\sc N_1+ \cdots + N_k = a \\
\sc n_1,\ldots,n_k \geq 0
\end{array} }
\frac{q^{N_1^2 + \cdots + N_k^2
+N_i+ \cdots + N_k}}
{(q)_{n_1} (q)_{n_2} \ldots (q)_{n_k}},
\end{equation}
which, up to a factor $(q)_a$, is the representation of
the Alder polynomials \cite{Alder} as found in ref.~\cite{Andrews74}.

\nsection{Proof of theorem~6}\label{pt4}
With the results of the previous two sections,
proving theorem~\ref{bt} is elementary.
First we define $S_{k,i,i';L}$ as the set of partitions
of $n$ of the form $n=\sum_{j=1}^{L-1} j f_j$ satisfying
the frequency conditions $f_1\leq i-1$,
$f_{L-1} \leq i'-1$ and
$f_j+f_{j+1}\leq k$ for $j=1,\ldots,L-2$.
Let $\pi$ be a partition in $S_{k,i,i';L}$,
with $\ell$ rows of length $L-1$.
Using the frequency condition this implies
$f_{L-2}\leq k-\ell$. Hence, by removing the
first $\ell$ rows, $\pi$ maps onto a partition
in $S_{k,i,k-\ell+1;L-1}$.
Conversely, by adding $\ell$ rows at the top
to a partition in $S_{k,i,k-\ell+1;L-1}$,
we obtain a partition in $S_{k,i,i';L}$.
Since in the above $\ell$ can take the values
$\ell=0,\ldots,i'-1$, the following recurrences
hold:
\begin{equation}
G_{k,i,i';L}(q) = \sum_{\ell=0}^{i'-1} q^{\ell(L-1)}
G_{k,i,k-\ell+1;L-1}(q) \qquad {\rm for} \quad i'=1,\ldots,k+1.
\label{Grec}
\end{equation}
In addition to this we have the initial condition
\begin{equation}
G_{k,i,i';2}(q) = \sum_{\ell=0}^{\min(i'-1,i-1)} q^{\ell}.
\end{equation}
Using the recurrence relations, it is in fact an easy matter
to verify that this is consistent with the condition
\begin{equation}
G_{k,i,i';0}(q) = \delta_{i,i'}.
\label{Ginit}
\end{equation}

Remains to verify that (\ref{Boson1}) and (\ref{Boson2}) satisfy the
recurrence (\ref{Grec}) and initial condition (\ref{Ginit}).
Since in these two equations we have used the variables $r$ and $s$
instead of $i'$ and $i$, let us first rewrite (\ref{Grec}) and
(\ref{Ginit}).
Suppressing the $k$, $s$ and $q$ dependence, setting
$G_{k,i,i';L}(q)=G_L(r)$, we get 
\begin{equation}
G_L(r) = \sum_{\ell=0}^{k-r} q^{\ell(L-1)}
G_{L-1}(\ell) \qquad {\rm for} \quad r=0,\ldots,k
\label{Grec2}
\end{equation}
and
\begin{equation}
G_0(r) = \left\{\begin{array}{lcl}
\delta_{s+r,k+1} \quad &{\rm for }&\;
s=1,3,\ldots,2\lfloor \frac{k}{2} \rfloor + 1 \\
\delta_{s-r,k+2} &{\rm for } & \;
s=2\lfloor \frac{k}{2} \rfloor + 3,\ldots,2k+1.
\end{array} \right.
\label{Ginit2}
\end{equation}

To verify that (\ref{Boson1}) and (\ref{Boson2}) satisfy the
initial condition (\ref{Ginit2}), we set $L=0$ and use the
fact that $r=0,1,\ldots,k$ and $s=1,3,\ldots,2k+1$.
From this and equation (\ref{Lis0})
one immediately sees that the only non-vanishing term in
(\ref{Boson1}) is given by $\Mults{0}{(k-s-r+1)/2}{r}{k}
=\delta_{s+r,k+1}$. Similarly the only non-vanishing term in
(\ref{Boson2}) is $\Mults{0}{(-k+s-r-2)/2}{r}{k}
=\delta_{s-r,k+2}$. 
Now recall that (\ref{Boson1}) with $L=0$ is 
$G_0(r)$ for $r\equiv k$.
From the allowed range of $r$ this implies 
$s=1,3,\ldots,2\lfloor \frac{k}{2} \rfloor + 1$, in
accordance with the top-line of (\ref{Ginit2}).
Also, since (\ref{Boson2}) with $L=0$ is
$G_0(r)$ for $r\not\equiv k$, and because of the range of
$r$, we get 
$s=2\lfloor \frac{k}{2} \rfloor + 3,\ldots,2k+1$, in accordance with the
second line in (\ref{Ginit2}).

Checking that (\ref{Boson1}) and (\ref{Boson2}) 
solve the recurrence relation (\ref{Grec2}) splits into
several cases due to the parity dependence of $G_L(r)$ and
of the $q$-multinomial recurrences (\ref{qrec}).
All of these cases are completely analogous and we
restrict our attention to $k$ and $r$ being even, so
that $G_L(r)$ is given by equation (\ref{Boson1}).
Substituting recurrences (\ref{qrec}), the first and
second sum in (\ref{qrec}) immediately give
the right-hand side of (\ref{Grec2}).
Consequently, the other two terms in (\ref{qrec}) give rise
to unwanted terms that have to cancel in order for (\ref{Grec2})
to be true.
Dividing out the common factor $q^{(2L-1)(k-r)/2}$  and
making the change of variables $m\to m-1$ in the last
sum of (\ref{qrec}), the unwanted terms read
\begin{eqnarray}
\lefteqn{
\renewcommand{\arraystretch}{0.6}
\sum_{\begin{array}{c}
\sc m=k-r+2 \\
\sc m \mbox{\scriptsize even}
\end{array} }^k
q^{-\frac{1}{2}m}
\sum_{j=-\infty}^{\infty} 
\left\{
q^{j\bigl((2j+1)(2k+3)-2s\bigr)}
\Mult{L-1}{\case{1}{2}(k L-s-m+1)+(2k+3)j}{m}{k} \right. }
\nonumber \\[-3mm]
&&\hspace{31mm}- 
q^{\bigl(2j+1\bigr)\bigl((2k+3)j+s\bigr)}
\Mult{L-1}{\case{1}{2}(k L+s-m+1)+(2k+3)j}{m}{k}
\nonumber \\
&&\hspace{31mm}+
q^{\bigl(2j-1\bigr)\bigl((2k+3)j-s\bigr)}
\Mult{L-1}{\case{1}{2}(k L+s-m+1)-(2k+3)j}{m-1}{k} 
\nonumber \\
&&\hspace{29mm}- \left.
q^{j\bigl((2j-1)(2k+3)+2s\bigr)}
\Mult{L-1}{\case{1}{2}(k L-s-m+1)-(2k+3)j}{m-1}{k} \right\}.
\end{eqnarray}
After changing the summation variable $j\to -j$ in the second and
fourth term, this becomes
\begin{eqnarray}
\lefteqn{
\renewcommand{\arraystretch}{0.6}
\sum_{\begin{array}{c}
\sc m=k-r+2 \\
\sc m \mbox{\scriptsize even}
\end{array} }^k
q^{-\frac{1}{2}m}
\sum_{j=-\infty}^{\infty} 
q^{j\bigl((2j+1)(2k+3)-2s\bigr)} \times} \\[-2mm]
\lefteqn{\hspace{-5mm}\left\{
\Mult{L-1}{\case{1}{2}(k L-s-m+1)+(2k+3)j}{m}{k}
- q^{s-2(2k+3)j}
\Mult{L-1}{\case{1}{2}(k L+s-m+1)-(2k+3)j}{m}{k} \right.}
\nonumber \\
\lefteqn{ \hspace{-5mm} \left.
+ q^{s-2(2k+3)j}
\Mult{L-1}{\case{1}{2}(k L+s-m+1)-(2k+3)j}{m-1}{k} 
-\Mult{L-1}{\case{1}{2}(k L-s-m+1)+(2k+3)j}{m-1}{k} \right\}.}
\nonumber 
\end{eqnarray}
We now show that the term within the curly braces vanishes
for all $m$ and $j$.
To establish this, we apply the symmetry (\ref{qs}) to all
four $q$-multinomials within the braces and divide by 
$q^{(k-m)(L-1)-\frac{1}{2}(kL-s-m+1)-(2k+3)j}$.
After replacing $L$ by $L+1$ and $m$ by $k-p$, this gives
\begin{eqnarray}
\lefteqn{
\Mult{L}{\case{1}{2}(k L+s-p-1)-(2k+3)j}{p}{k}
-\Mult{L}{\case{1}{2}(k L-s-p-1)+(2k+3)j}{p}{k}}
\nonumber \\
\lefteqn{+ q^L
\Mult{L}{\case{1}{2}(k L-s-p-1)+(2k+3)j}{p+1}{k} 
-q^L
\Mult{L}{\case{1}{2}(k L+s-p-1)-(2k+3)j}{p+1}{k} .}
\end{eqnarray}
Recalling the tautology (\ref{taut}) with $a=\frac{1}{2}(kL+s-p-1)
-(2k+3)j$ this indeed gives zero.

\nsection{Proof of theorem~5}\label{pt3}
\subsection{From partitions to paths}
To prove expression (\ref{Fermion}) of theorem~\ref{ft}, 
we reformulate the problem of calculating the generating
function $G_{k,i,i';L}(q)$ into a lattice path
problem. Hereto we represent each partition $\pi$
in $S_{k,i,i';L}$ as a restricted lattice path $p(\pi)$,
similar in spirit to the lattice path formulation of
the left-hand side of (\ref{An}) by 
Bressoud~\cite{Bressoud}.\footnote{Finitizing 
Bressoud's lattice paths by fixing
the length of his paths to $L$,
results in the left-hand side
of (\ref{fq}) of the next section. 
Hence the lattice paths introduced here are intrinsically
different from those of ref.~\cite{Bressoud} and in fact
correspond to a finitization of the paths of ref.~\cite{RV}.}

\begin{figure}[t]
\epsfxsize = 13cm
\centerline{\epsffile{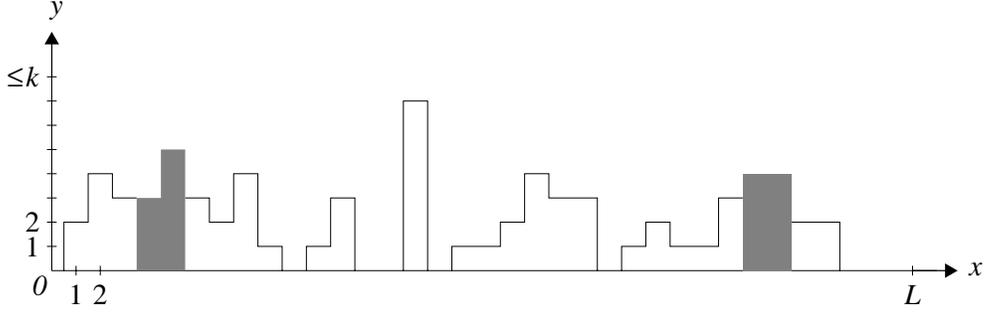}}
\caption{A lattice path of the partition $(f_1,\ldots,f_{L-1})=
(2,4,3,3,5,3,2,4,1,0,1,3,0,0,7,0,$
$1,1,2,4,3,3,0,1,2,1,1,3,4,4,2,2,0,0)$. 
The shaded regions correspond to the two
particles with largest charge (=8), as described below.}
\label{latpath}
\end{figure}

To map a partition $\pi$ of $n=\sum_{j=1}^{L-1} j f_j$
onto a lattice path $p(\pi)$, draw a horizontal 
line-segment in the $(x,y)$-plane
from $(j-\frac{1}{2},f_j)$ to $(j+\frac{1}{2},f_j)$ 
for each $j=1,\ldots,L-1$. 
Also draw vertical line-segments  from $(j+\frac{1}{2},f_j)$ to
$(j+\frac{1}{2},f_{j+1})$ for all $j=0,\ldots,L-1$, 
where $f_0=f_{L-1}=0$.
As a result $\pi$ is represented by a lattice path (or
histogram) from $(\case{1}{2},0)$ to $(L-\case{1}{2},0)$.
The frequency condition $f_j+f_{j+1}\leq k$
translates into the condition that the sum of the {\em heights}
of a path at $x$-positions $j$ and $j+1$ does not exceed $k$.
The restrictions $f_1 \leq i-1$ and $f_{L-1} \leq i'-1$
correspond to the restrictions that the heights at $x=1$ and $x=L-1$
are less than $i$ and $i'$, respectively.
An example of a lattice path for $k\geq 8$, $i\geq 3$ and
$i'\geq 1$, is shown in Figure~\ref{latpath}.

The above map clearly is reversible, and any lattice path
satisfying the above height conditions maps onto a partition
in $S_{k,i,i';L}$. From now on we let $P_{k,i,i';L}$ denote
the set of restricted lattice paths corresponding to 
the set of partions $S_{k,i,i';L}$.

From the map of partitions onto paths, the
problem of calculating the generating function 
$G_{k,i,i';L}(q)$ can be reformulated as
\begin{equation}
G_{k,i,i';L}(q) = 
\sum_{\textstyle p\!\in\!\! P_{k,i,i';L}} W(p) 
\label{Gpath}
\end{equation}
with Boltzmann weight $W(p)=\prod_{j=1}^{L-1} q^{j f_j}$.

Before we actually compute the above sum, we remark 
that in the following $k,i$ and $i'$ will always be fixed.
Hence, to simplify notation, we use $G_L$ and $P_L$ to denote
$G_{k,i,i';L}$ and $P_{k,i,i';L}$, respectively.

\subsection{Fermi-gas partition function; $i=i'=k+1$}
To perform the sum (\ref{Gpath}) over the
restricted lattice path, we follow a procedure similar to
the one employed in our proof of Virasoro-character identities
for the unitary minimal models~\cite{Wa,Wb}.
That is, the sum (\ref{Gpath}) is interpreted as the
grand-canonical partition function of a one-dimensional
lattice-gas of fermionic particles.

The idea of this approach is to view each lattice path
as a configuration of particles on a one-dimensional lattice.
Since not all lattice paths correspond to the same
{\em particle content} $\vec{n}$, this gives rise to a natural
decomposition of (\ref{Gpath}) into
\begin{equation}
G_L(q)  = \sum_{\vec{n}} Z_L(\vec{n};q),
\label{GZ}
\end{equation}
with $Z_L$ the canonical partition function,
\begin{equation}
Z_L(\vec{n};q) = 
\sum_{\textstyle p\!\in\!\! P_L(\vec{n})} W(p).
\end{equation}
Here $P_L(\vec{n})\subset P_L$ is the set of paths corresponding to a
particle configuration with content $\vec{n}$.
To avoid making the following description of the lattice gas unnecessarily
complicated, we assume $i=i'=k+1$ in the remainder of this section.
Subsequently we will briefly indicate how to modify the calculations
to give results for general $i$ and $i'$.

\begin{figure}[t]
\epsfxsize = 15cm
\centerline{\epsffile{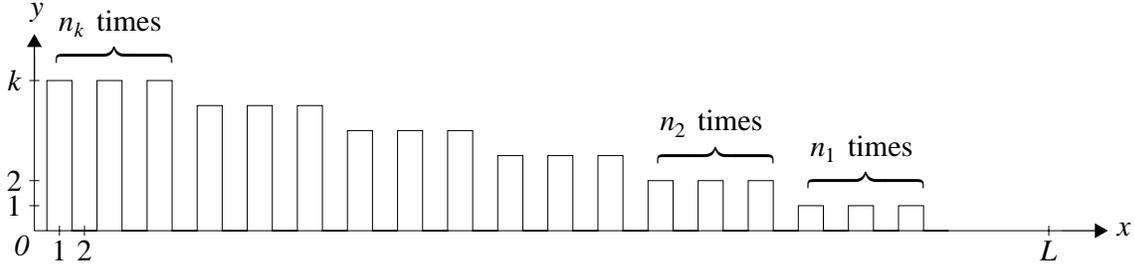}}
\caption{The minimal path of content $\vec{n}=(n_1,\ldots,n_k)^T$.}
\label{min}
\end{figure}

To describe how to interpret each path in $P_L=P_{k,k+1,k+1;L}$
as a particle configuration,
we first introduce a special kind of paths from which all 
other paths can be constructed.
\begin{definition}[minimal paths]
The path shown in Figure~\ref{min} is called 
the minimal path of content $\vec{n}$.
\end{definition}
\begin{definition}[charged particle]
In a minimal path, each column with non-zero height $t$ corresponds to a 
particle of charge $t$.
\end{definition}
Note that in the minimal path the particles are ordered according to
their charge and that adjacent particles are separated by
a single empty column.
The number of particles of charge $t$ is denoted $n_t$ and $\vec{n}=
(n_1,\ldots,n_k)^T$.
For later use it will be convenient to give each particle a label,
$p_{t,\ell}$ denoting the $\ell$-th particle of charge $t$,
counted from the right.

Since the length of a path is fixed by $L$,
there are only a finite number of minimal paths.
In particular, we have $2(n_1+\cdots + n_k) \leq L$, so that there
are $
\mults{\lf L/2 \rf +k}{k}
$
different minimal paths.

In the following we show that all non-minimal paths in $P_L$ can 
be constructed out of one (and only one) minimal path using a set 
of elementary moves.
Hereto we first describe how various local configurations may be 
changed by moving a particle.\footnote{In moving a particle we always
mean motion from left to right.}
To suit the eye, the particle being moved in each example has been shaded.

To describe the moves we first consider the simplest type of
motion, when the two columns immediately to the right of a particle 
are empty.
\begin{definition}[free motion]
The following sequence of moves is called free motion:
\[ \hspace{3cm} \epsfxsize = 9cm
\epsffile{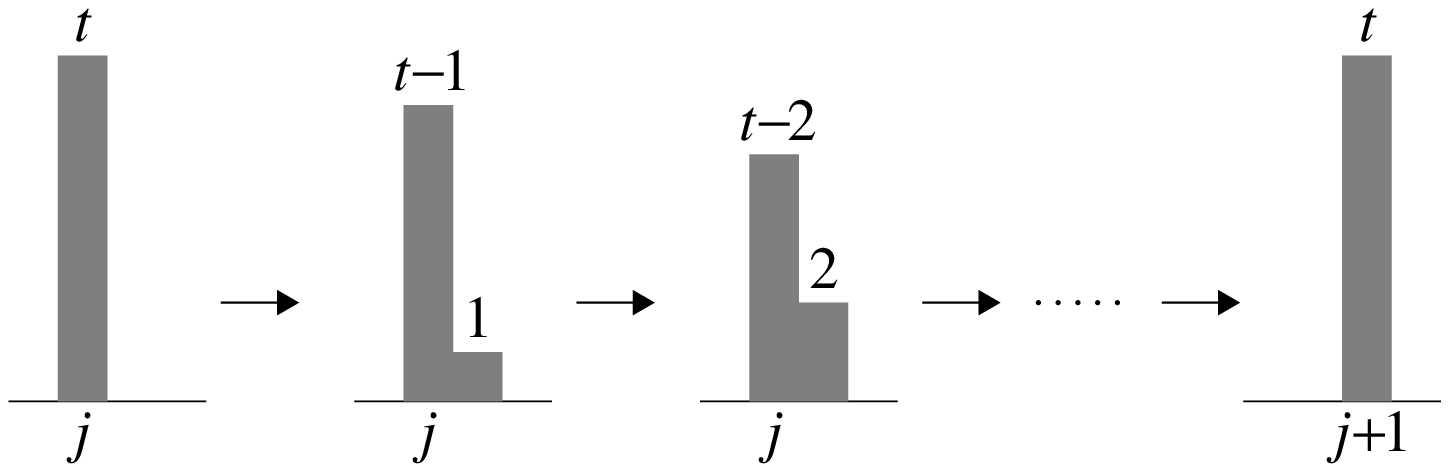} \]
\end{definition}
Clearly, a particle of charge $t$ in free motion takes 
$t$ moves to fully shift position by one unit.

\begin{figure}[h]
\epsfxsize = 12cm
\centerline{\epsffile{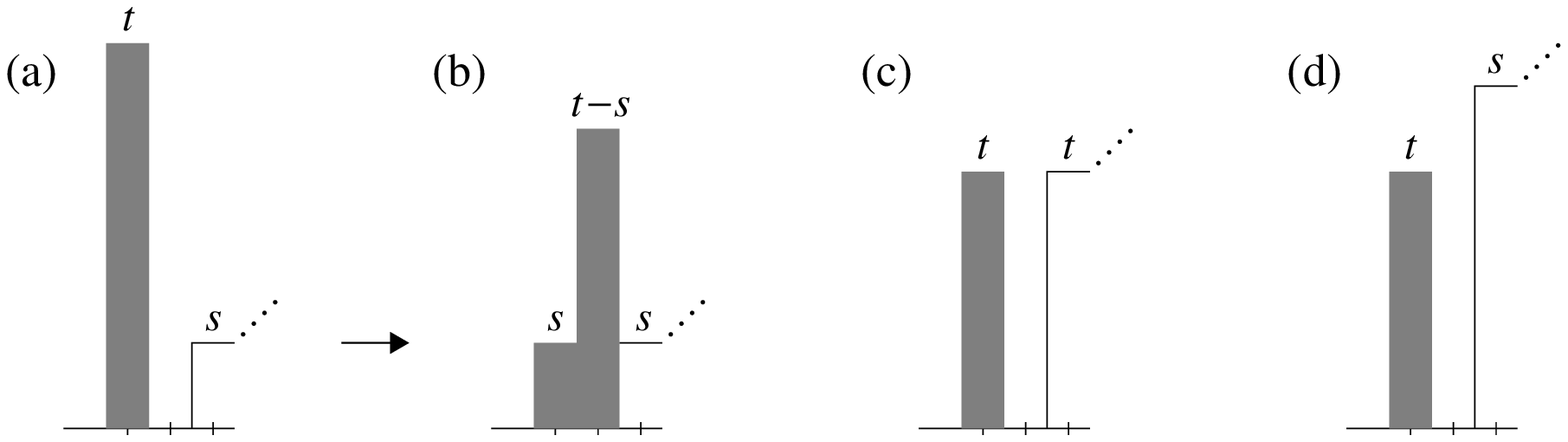}}
\caption{}
\label{ts}
\end{figure}

Now assume that in moving a particle of charge $t$,
we at some stage encounter the local configuration shown in 
Figure~\ref{ts}(a). We then allow the particle to make $t-s$ more
moves following the rules of free motion, to end up with
the local configuration shown in (b).
If instead of (a) we encounter the configurations (c) or (d),
the particle can make no further moves.

\begin{figure}[h]
\epsfxsize = 15cm
\centerline{\epsffile{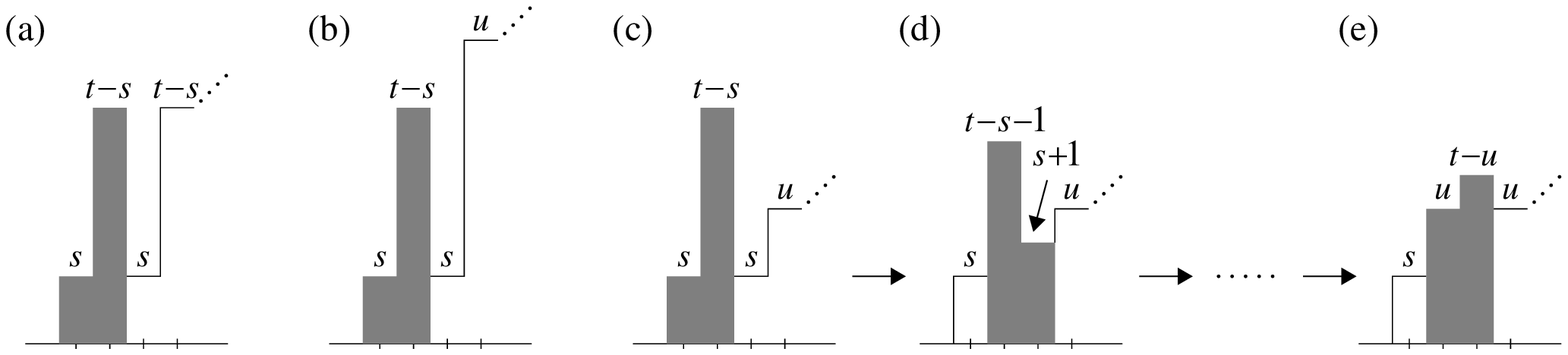}}
\caption{}
\label{ts2}
\end{figure}

In case of the configuration of Figure~\ref{ts}(b), there are
three possibilities. Either we have one of 
the configurations shown in Figure~\ref{ts2}(a) and (b),
in which case the particle cannot move any further,
or we have the configuration shown in Figure~\ref{ts2}(c) 
(with $0\leq u < t-s$), 
in which case we can make $t-u-s$ moves, going from
(c) to (e).
Ignoring the for our rules irrelevant column 
immediately to the left of the
particle, the configuration of Figure~\ref{ts2}(e)
is essentially the same as that of Figure~\ref{ts}(b).
To further move the particle we can thus refer to the rules
given in Figure~\ref{ts2}. That is, if the column immediately
to the right of the white column of height $u$ has height $t-u$ 
(corresponding to the configuration of Figure~\ref{ts2}(a) with $s\to u$)
the particle cannot move any further.
Similarly, if the column immediately
to the right of the white column of height $u$ has height $v>t-u$ 
(corresponding to the configuration of 
Figure~\ref{ts2}(b) with $u\to v$ and $s\to u$)
the particle cannot move any further.
However, if the height of the column immediately to the right
of the white column of height $u$ is $0\leq v<t-u$
(corresponding to the configuration of Figure~\ref{ts2}(c) 
with $u\to v$ and $s\to u$) we can make another $t-u-v$ moves.

Having introduced all necessary moves we come to the main propositions
of this section
\begin{proposition}[rules of motion]\label{prm}
Each non-minimal configuration can be obtained from one and only one
minimal configuration by letting the particles carry out elementary moves
in the following order:
\begin{itemize}
\item
Particle $p_{t,\ell}$ moves prior to $p_{s,\ell'}$ if $t<s$.
\item
Particle $p_{t,\ell'}$ moves prior to $p_{t,\ell}$ if $\ell'<\ell$.
\end{itemize}
\end{proposition}
To prove this let us assume that we have completed the motion
of all particles of charge less than $t$ and all particles $p_{t,\ell'}$ 
with $\ell'<\ell$, and that we are
currently moving the particle $p_{t,\ell}$.
Therefore, for $x\leq 2n_k+\cdots + 2n_{t+1} + 2(n_t-\ell):=x_{\min}$,
the lattice path still corresponds to the minimal path.

Now note that in moving the particle $p_{t,\ell}$,
we never create a local configuration in which the sum of the height
of two consecutive columns is greater than $t$, see
the free motion and the Figures~\ref{ts}(a)(b) and \ref{ts2}(c)--(e).
Since we have not moved any of the particles of charge greater than
$t$, this means that to the right of $x_{min}$ no two
consecutive columns have summed heights greater than $t$.\footnote{
This also means that in moving a particle of charge
$t$, the configurations of Figure~\ref{ts}(d) and Figure~\ref{ts2}(b)
in fact never arise.}
Also note that as soon as $p_{t,\ell}$ meets two columns immediately
to its right whose summed heights equal $t$, $p_{t,\ell}$ cannot
move any further, see Figures~\ref{ts}(c) and ~\ref{ts2}(a).
Consequently, $p_{t,\ell}$ always corresponds to the
leftmost two consecutive columns right of $x_{\min}$
whose summed height equals $t$. (In fact, it corresponds to
the leftmost two consecutive columns right of $x_{min}$ with maximal
summed heights.)

Now we define reversed moves by reading all the previous 
figures in a mirror. Using this motion we can move
$p_{t,\ell}$ all the way back to its minimal position but
not any further.
To see this we note that the only situations in which $p_{t,\ell}$ 
cannot be moved further back is if it meets two 
consecutive columns to its left whose summed heights
are greater or equal to $t$. Since we have just argued that such
a configuration cannot occur between $x_{\min}$ and $p_{t,\ell}$
we can indeed move $p_{t,\ell}$ back to its minimal position using
the reversed moves. Once it is back in its minimal position
we either have the mirror image of Figure~\ref{ts}(c) (in case $\ell<n_t$)
or (d) (in case $\ell=n_t$). Neither of these configurations allows
for further reversed moves. 

The above, however, gives a general procedure for reducing each non-minimal
path to a minimal one by simply reversing the rules of motion
in the proposition.
That is, we first scan the path for all particles of charge $k$,
by locating all occurrence of two consecutive columns of
summed heights $k$. From left to right these label the particles $p_{k,n_k}$
to $p_{k,1}$. Applying the previous reasoning with $t=k$, we can
first move $p_{k,n_k}$ back using reversed moves, than  
$p_{k,n_k-1}$, et cetera, until  all particles of charge $k$
have taken their ``minimal position''.
Repeating this for the particles of charge $k-1$, then
the particles of charge $k-2$, et cetera, each non-minimal
path reduces to a unique minimal path.

\begin{figure}[t]
\epsfxsize = 13cm
\centerline{\epsffile{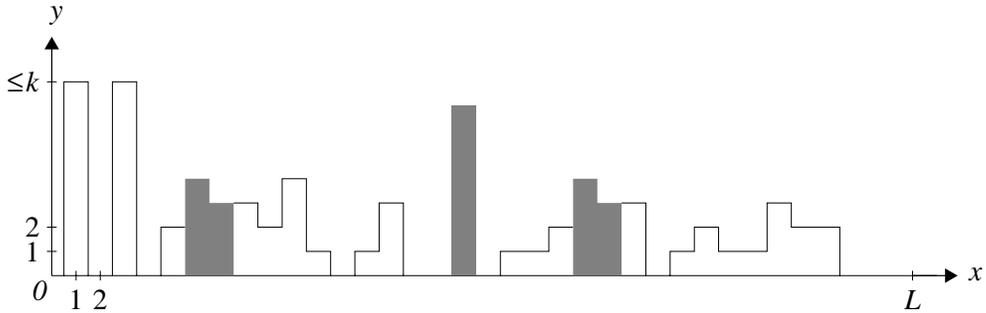}}
\caption{The lattice path obtained from Figure~1 after
moving the largest particles to their minimal position.
The shaded regions mark the three next-largest particles.}
\label{redpath}
\end{figure}

As an example to the above,
for the path of Figure~\ref{latpath} the shaded regions
mark the (two) particles with largest charge (=8). Moving them
back using the reversed motion, the leftmost particle being
moved first, we end up with the path shown in Figure~\ref{redpath},
in which now the (three) particles with next-largest charge 
have been marked. We leave it to the reader to further reduce
the path to obtain the minimal path of content $(n_1,n_2,\ldots,n_8)=
(2,1,2,1,3,1,3,2)$.

The elementary moves and the reversed moves are clearly reversible.
If a particle of charge $t$ has made an elementary move changing a path
from $p$ to $p'$, we can always carry a reversed move going 
from $p'$ back to $p$.
Since each path can be reduced to a unique minimal path using the
reversed moves by carrying out the rules of motion of proposition~\ref{prm}
in reversed order, we have thus established that using the rules of motion
we can generate each non-minimal path uniquely from a minimal path.
Hence the proposition is proven.

We now have established the decomposition of the sum (\ref{Gpath}) into
(\ref{GZ}), where $Z_L(\vec{n};q)$ is the generating function of the
paths generated by the minimal path labelled by $\vec{n}$,
or, in other words, $Z_L(\vec{n};q)$ is the partition function of
a lattice gas of fermions with particle content $\vec{n}$.
The fermionic nature being that, unlike particles of different charge,
particles of equal charge cannot exchange position.

Our next result concerns the actual computation of  the
partition function.
\begin{proposition}\label{pZ}
The partition function $Z_L$ is given by
\begin{equation}
Z_L(\vec{n};q) = 
q^{\displaystyle \, \vec{n}^T C^{-1}_k \vec{n}} 
\prod_{t=1}^k
\Bin{n_t+m_t}{n_t},
\label{ZLn}
\end{equation}
with {\rm $
\vec{m}+\vec{n}=
\case{1}{2}({\cal I}_k \, \vec{m}
+L \, \vec{\e}_k)$}.
\end{proposition}
To prove this we first determine
the contribution to $Z_L$ of the minimal path of content $\vec{n}$,
\begin{eqnarray}
W(p_{\min})/\ln q  &=& 
\sum_{t=1}^k t \sum_{\ell=1}^{n_t} 
\Bigl(2\ell-1 + 2 \sum_{s=t+1}^k n_s\Bigr)
=
\sum_{t=1}^k t n_t
\Bigl(n_t + 2 \sum_{s=t+1}^k n_s\Bigr) \nonumber \\
&=& 
\sum_{r=1}^k \sum_{t=r}^k n_t 
\Bigl(n_t + 2 \sum_{s=t+1}^k n_s\Bigr) 
= 
\sum_{r=1}^k N_r^2 = 
\vec{n}^T C^{-1}_k \vec{n} .
\label{Wmin}
\end{eqnarray}

To obtain the contribution to $Z_L$ of the non-minimal
configurations, we apply the rules of motion of proposition~\ref{prm}.
If $e_{\ell}$ denotes the number of elementary moves carried out by
$p_{t,\ell}$, the generating function of moving the particles
of charge $t$ reads
\begin{equation}
\sum_{e_1=0}^{m_t} 
\sum_{e_2=0}^{e_1}
\ldots
\sum_{e_{n_t}=0}^{e_{n_t-1}}
q^{e_1+e_2 + \cdots + e_{n_t}} = \Bin{m_t + n_t}{n_t}.
\label{respart}
\end{equation}
Here we have used the fact that each elementary move generates a factor $q$
and that $p_{t,\ell}$ cannot carry out more elementary moves than
$p_{t,\ell-1}$. (If $p_{t,\ell}$ has made as many moves as $p_{t,\ell-1}$
we obtain either the local configuration of Figure~\ref{ts}(c) or
Figure~\ref{ts2}{a}, prohibiting any further moves.)
The number $m_t$ in (\ref{respart}) is the maximal number of
elementary moves $p_{t,1}$ can make and remains to be determined.
If the content of the minimal path is $\vec{n}$,
the $x$-position of $p_{t,1}$ in $p_{\min}$ is $2(n_k+\cdots + n_t)-1:=x_0$.
To fix $m_t$, let us assume that after the motion of the particles 
of charge less than $t$ has been completed, the nontrivial part of the 
lattice path is encoded by the sequence of heights
$(f_{x_0},\ldots,f_{L-1})$, with $f_{x_0}=t$, $f_{x_0+1}=0$ and
$f_j+f_{j+1}<t$ for $j>x_0$.
The particle $p_{t,1}$ can now move all the way to $x=L-1$ 
making
\begin{eqnarray}
m_t &=&
(t-f_{x_0+2}) +
(t-f_{x_0+2}-f_{x_0+3}) +
(t-f_{x_0+3}-f_{x_0+4}) \nonumber \\ 
& & \qquad \qquad + \cdots
+ (t-f_{L-2}-f_{L-1}) + (t-f_{L-1}) \nonumber \\
&=& t(L-x_0-1) -2\sum_{j=x_0+2}^{L-1} f_j
\label{mt}
\end{eqnarray}
elementary moves.
To simplify this, note that the sum on the right-hand side is
nothing but twice the sum of the heights of the columns
right of $x=x_0$, which is $2 \sum_{s=1}^{t-1} s n_s$.
Substituting this in (\ref{mt}) and using the definition of
$x_0$, results in 
\begin{equation}
m_t  = t L - 2\sum_{s=1}^{t-1} s n_s - 2 t \sum_{s=t}^k n_t
= t L -2\sum_{s=1}^k \min(s,t) n_s
= L(C^{-1}_k)_{t,k}
-2\sum_{s=1}^k (C^{-1}_k)_{t,s} n_s
\label{mt2}
\end{equation}
in accordance with proposition~\ref{pZ}.

Putting together the results (\ref{Wmin}), (\ref{respart})
and (\ref{mt2}) completes the proof of proposition~\ref{pZ}.
Substituting the form (\ref{ZLn}) of the partition function 
into (\ref{GZ}) proves expression
(\ref{Fermion}) of theorem~\ref{ft}, for $i=i'=k+1$.

\subsection{Fermi-gas partition function; general $i$ and $i'$}
Modifying the proof of theorem~\ref{ft} for $i=i'=k+1$ to all
$i$ and $i'$ is straightforward and few details will be given.
It is in fact interesting to note that unlike our proof for the
character identities of the unitary minimal models~\cite{Wa,Wb}, 
the general case here does not require the introduction of 
additional ``boundary particles''. 

First let us consider the general $i'$ case, with $i=k+1$.
This implies that the height $f_{L-1}$ of the last column
of the lattice paths is no longer free to take any of the values
$1,\ldots,k$, but is bound by $f_{L-1}\leq i'-1$.
For the particles of charge less or equal to $i'-1$ this does
not impose any new restrictions on the maximal number of moves
$p_{t,1}$ can make. For $t>i'-1$ however, $m_t$ in (\ref{mt2})
has to be decreased by $t-i'+1$.
Thus we find that $m_t$ of (\ref{mt2}) has to be replaced by
$m_t - \max(0,t-i'+1)= m_t -t + \min(t,i'-1)$.
Recalling $(C^{-1}_k)_{s,t}=\min(s,t)$, this yields
\begin{equation}
m_t = (L-1)(C^{-1}_k)_{t,k} 
+ (C^{-1}_k)_{t,i'-1} 
-2\sum_{s=1}^k (C^{-1}_k)_{t,s} n_s
\end{equation}
and therefore,  $m_t+n_t = 
\case{1}{2}(\sum_{s=1}^k ({\cal I}_k)_{t,s}  m_s
+(L-1) \delta_{t,k} +\delta_{t,i'-1})$
which is in accordance with proposition~\ref{pZ}, for $i=k+1$.

\begin{figure}[t]
\epsfxsize = 15cm
\centerline{\epsffile{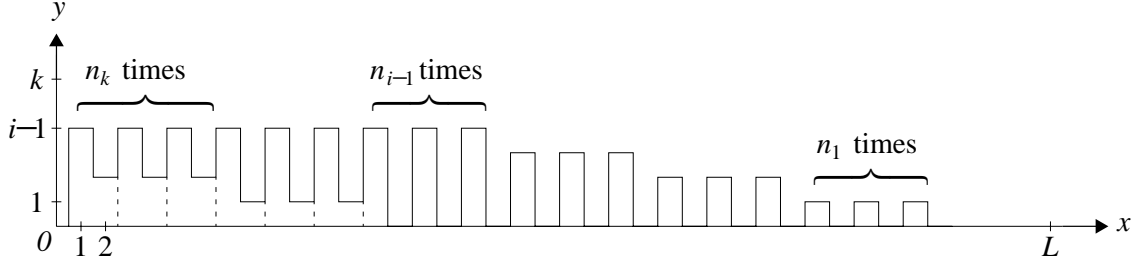}}
\caption{The minimal path of content $\vec{n}=(n_1,\ldots,n_k)^T$
for general $i$. The dashed lines are drawn to mark the
different particles.}
\label{min2}
\end{figure}

Second, consider $i$ general, but $i'=k+1$, so that
$f_1\leq i-1$, $f_{L-1}\leq k$.
Now the modification is slightly more involved
since the actual minimal paths change from those of Figure~\ref{min}
to those of Figure~\ref{min2}.
This leads to a change in the calculation of $W(p_{\min})$ to
\begin{eqnarray}
W(p_{\min})/\ln q  &=& 
\sum_{t=1}^k \sum_{\ell=1}^{n_t} 
\Bigl(2\ell t - \min(t,i-1) + 2t \sum_{s=t+1}^k n_s\Bigr)
\nonumber \\
&=&
\sum_{t=1}^k t n_t
\Bigl(n_t + 2 \sum_{s=t+1}^k n_s\Bigr)
+\sum_{t=i}^k (t-i+1) n_t \nonumber \\
&=& 
\vec{n}^T C^{-1}_k \vec{n} + \sum_{t=1}^k t n_t -
\sum_{t=1}^k \min(t,i-1) n_t \; = \; 
\vec{n}^T C^{-1}_k (\vec{n} + \e_k -\e_{i-1}),
\end{eqnarray}
which is indeed the general form of the quadratic exponent
of $q$ in (\ref{Fermion}).
Also $m_t$ again requires modification, which is in fact
similar to the previous case:
$m_t \to m_t - \max(0,t-i+1)$.
To see this note that it takes $\max(0,t-i+1)$
elementary moves to move $p_{t,1}$ from its
minimal position in Figure~\ref{min} to its minimal
position in Figure~\ref{min2}.

Finally, combining the previous two cases, and using the fact
that the modifications of $m_t$ due to $f_1\leq i-1$ and 
$f_{L-1}\leq i'-1$ 
are independent, we immediately arrive at the general form of
(\ref{Fermion}) with $(m,n)$-system (\ref{mn}).

\nsection{Discussion}\label{secdis}
In this paper we have presented polynomial identities 
which arise from finitizing Gordon's frequency partitions.
The bosonic side of the identities involves 
$q$-deformations of the coefficient of $x^a$ in the
expansion of $(1+x+\cdots +x^k)^L$.
The fermionic side follows from interpreting
the generating function of the frequency partitions,
as the grand-canonical partition function of a one-dimensional
lattice gas.

Interestingly, in recent publications, Foda and Quano, and Kirillov,
have given different polynomial identities which imply 
(\ref{An})~[17-19].
In the notation of section~\ref{intro} these identities 
can be expressed as
\begin{theorem}[Foda--Quano, Kirillov]
For all $k\geq 1$, $1\leq i\leq k+1$ and $L\geq k-i+1$,
\begin{eqnarray}
\lefteqn{
\sum_{n_1,\ldots,n_k\geq 0}
q^{N_1^2 + \cdots + N_k^2 + N_i + \cdots + N_k}
\prod_{j=1}^{k}
\Bin{L-N_j-N_{j+1}-2\sum_{\ell=1}^{j-1} N_{\ell}
-\alpha_{i,j}}{n_j} } \nonumber \\
&&\qquad = 
\sum_{j=-\infty}^{\infty} (-)^j 
q^{j\bigl( (2k+3)(j+1)-2i\bigr)/2}
\Bin{L}{\lfloor \frac{L-k+i-1-(2k+3)j}{2} \rfloor},
\label{fq}
\end{eqnarray}
with $N_{k+1}=0$ and $\alpha_{i,j}=\max(0,j-i+1)$.
\end{theorem}

An explanation for this different finitization of (\ref{An})
can be found in a theorem due to Andrews~\cite{Andrews72}:
\begin{theorem}[Andrews]
Let $Q_{k,i}(n)$ be the number of partitions of $n$
whose successive ranks lie in the interval
$[2-i,2k-i+1]$.
Then $Q_{k,i}(n)=A_{k,i}(n)$.
\end{theorem}
It turns out that it is the (natural) finitization of these
successive rank partitions which gives rise to the
above alternative polynomial finitization. That is,
(\ref{fq}) is an identity for the generating function
of partitions with largest part $\leq 
\lfloor (L+k-i+2)/2 \rfloor$, number of parts
$\leq \lfloor (L-k+i-1)/2 \rfloor$, whose successive ranks 
lie in the interval $[2-i,2k-i+1]$.

Let us now reexpress (\ref{fq}) into a form similar
to equations (\ref{Fermion}), (\ref{Boson1}) and (\ref{Boson2}). 
Hereto we eliminate $i$ in the right-hand side of
(\ref{fq}) in favour of the variable $s$ of equation
(\ref{defs}) and split the result into two cases.
This gives
\begin{eqnarray}
\lefteqn{\mbox{RHS}(\ref{fq})=
\sum_{j=-\infty}^{\infty} \left\{
q^{j\bigl((2j+1)(2k+3)-2s\bigr)}
\Mult{L}{\case{1}{2}(L+k-s+1)+(2k+3)j}{0}{1}
\right. } \nonumber \\
& &  \qquad \qquad \qquad \left.
-q^{\bigl(2j+1\bigr)\bigl((2k+3)j+s\bigr)}
\Mult{L}{\case{1}{2}(L+k+s+1)+(2k+3)j}{0}{1}
\right\}
\end{eqnarray}
for $L+k$ even, and
\begin{eqnarray}
\lefteqn{\mbox{RHS}(\ref{fq})=
\sum_{j=-\infty}^{\infty} \left\{
q^{j\bigl((2j+1)(2k+3)-2s\bigr)}
\Mult{L}{\case{1}{2}(L-k+s-2)-(2k+3)j}{0}{1}
\right. } \nonumber \\
& &  \qquad \qquad \qquad \left.
-q^{\bigl(2j+1\bigr)\bigl((2k+3)j+s\bigr)}
\Mult{L}{\case{1}{2}(L-k-s-2)-(2k+3)j}{0}{1}
\right\}
\end{eqnarray}
for $L+k$ odd.
This we recognize to be exactly (\ref{Boson1}) and (\ref{Boson2}) with
$r=0$, $k L$ replaced by $L$ and $\Mults{\ldots}{\ldots}{0}{k}$
replaced by $\Mults{\ldots}{\ldots}{0}{1}$.

Similarly, if we express the left-hand side of (\ref{fq})
through an $(n,m)$-system, we find 
precisely (\ref{Fermion}) but with
\begin{equation}
\vec{m}+\vec{n}=
\case{1}{2}\Bigl({\cal I}_k \, \vec{m}
+L \, \vec{\e}_1 + \vec{\e}_{i-1} -\vec{\e}_k \Bigr).
\label{mnFQ}
\end{equation}
This is just (\ref{mn}) with $r=0$ and 
$L\, \e_k$ replaced by $L\, \e_1$.

From the above observations it does not require much insight to
propose more general polynomial identities which have
(\ref{fq}) and those implied by the theorems~\ref{ft} and
\ref{bt} as special cases.
In particular, we have confirmed the following conjecture 
by extensive series expansions.
\begin{conjecture}\label{conj}
For all $k\geq 1$, $1\leq \ell \leq k$, $1\leq i \leq k+1$,
$1\leq i'\leq \ell+1$ and $\ell L\geq k+\ell-i-i'+2$
{\rm
\begin{equation}
\sum_{n_1,n_2,\ldots,n_k\geq 0}
q^{\displaystyle \, \vec{n}^T C^{-1}_k (\vec{n}
+\vec{\e}_k -\vec{\e}_{i-1})}
\prod_{j=1}^k
\Bin{n_j+m_j}{n_j}
\end{equation}}
with $(m,n)$-system given by
{\rm
\begin{equation}
\vec{m}+\vec{n}=
\case{1}{2}\Bigl({\cal I}_k \, \vec{m}
+(L\!-1) \, \vec{\e}_{\ell} + \vec{\e}_{i-1}
+\vec{\e}_{i'-1}-\vec{\e}_k \Bigr)
\label{mnGen}
\end{equation}}
equals
\begin{eqnarray}
\lefteqn{
\sum_{j=-\infty}^{\infty} \left\{
q^{j\bigl((2j+1)(2k+3)-2s\bigr)}
\Mult{L}{\case{1}{2}(\ell L+k-s-r+1)+(2k+3)j}{r}{\ell} 
\right. } \nonumber \\
& & \hspace{2mm} \left.
-q^{\bigl(2j+1\bigr)\bigl((2k+3)j+s\bigr)}
\Mult{L}{\case{1}{2}(\ell L+k+s-r+1)+(2k+3)j}{r}{\ell}
\right\}
\end{eqnarray}
for $r\equiv \ell L+k \mod{2}$ and
\begin{eqnarray}
\lefteqn{
\sum_{j=-\infty}^{\infty} \left\{
q^{j\bigl((2j+1)(2k+3)-2s\bigr)}
\Mult{L}{\case{1}{2}(\ell L-k+s-r-2)-(2k+3)j}{r}{\ell} 
\right. } \nonumber \\
& & \hspace{2mm} \left.
-q^{\bigl(2j+1\bigr)\bigl((2k+3)j+s\bigr)}
\Mult{L}{\case{1}{2}(\ell L-k-s-r-2)-(2k+3)j}{r}{\ell}
\right\}
\end{eqnarray}
for $r\not\equiv \ell L+k \mod{2}$.
Here $s$ is defined as in (\ref{defs}) and 
\begin{equation}
r=\ell-i'+1
\end{equation}
so that $r=0,\ldots,\ell$.
\end{conjecture}
For later reference, let us denote these more general polynomials
as $G_{k,i,i';L}^{(\ell)}(q)$. Then $\ell=k$ corresponds
to the polynomials considered in this paper and
$\ell=1$ to those of Foda, Quano and Kirillov.

The above conjecture leads one to wonder whether there are in fact
(at least) $k$ different partition theoretical interpretations 
of (\ref{An}), each of which has a natural finitization
corresponding to the polynomials $G_{k,i,i';L}^{(\ell)}(q)$
with $\ell=1,\ldots,k$.

Intimately related to the conjecture and perhaps even more surprising
is the following observation, originating from the work
of Andrews and Baxter~\cite{AB}.
For $k\geq 0$ and $1\leq i \leq k+1$,
define a $k$-variable generating function
\begin{equation}
f(x_1,\ldots,x_k) = 
\sum_{n_1,n_2,\ldots,n_{k}\geq 0}
\frac{q^{N_1^2 + \cdots + N_k^2+N_i+\cdots+N_k} 
\, x_1^{2N_1} x_2^{2(N_1+N_2)} 
\cdots x_k^{2(N_1+\cdots + N_k)}}
{(x_1)_{n_1+1} (x_2)_{n_2+1} \cdots (x_k)_{n_k+1} } \: ,
\label{fk}
\end{equation}
where $(x)_n = \prod_{k=0}^{n-1} (1-xq^k)$.
Obviously, $(1-x_1) \cdots (1-x_k) f(1,\ldots,1)$
corresponds to the left-hand side of (\ref{An}). 
Now define the polynomials $P(\ell_1,\ldots,\ell_k):=P(\vec{\ell})$
as the coefficients in the series expansion of $f$,
\begin{equation}
f(x_1,\ldots,x_k) = \sum_{\ell_1,\ldots,\ell_k} 
P(\vec{\ell}) \,x_1^{\ell_1} \cdots x_k^{\ell_k}.
\end{equation}
From the readily derived functional equations for $f$
and the recurrences (\ref{Grec}) with (\ref{Ginit})
one can deduce that
\begin{equation}
P(\vec{m}+2 C_k^{-1} \vec{n}) = G_{k,i,i';L}(q)
\end{equation}
with $\vec{m}$ and $\vec{n}$ given by (\ref{mn}).
Similarly the polynomials of Foda, Quano and Kirillov arise
again as $P(\vec{m}+2 C_k^{-1} \vec{n})$ where $\vec{m}$ and
$\vec{n}$ now satisfy (\ref{mnFQ}).
Again we found numerically that also the polynomials featuring the
conjecture appear naturally. That is,
\begin{equation}
P(\vec{m}+2 C_k^{-1} \vec{n}) = G^{(\ell)}_{k,i,i';L}(q),
\end{equation}
where now the generalized $(m,n)$-system (\ref{mnGen})
should hold (so that $\vec{m}+2 C_k^{-1} \vec{n}=
C^{-1}_k ((L\!-1) \, \vec{\e}_{\ell} + \vec{\e}_{i-1}
+\vec{\e}_{i'-1}-\vec{\e}_k)$).

\vspace*{5mm}

Although all the polynomial identities implied by 
conjecture~\ref{conj} reduce to Andrews' identity~(\ref{An})
in the $L$ to infinity limit,
they still provide a powerful tool for generating
new $q$-series results.
That is, if we first replace $q\to 1/q$ and then take $L\to\infty$,
new identities arise.
To state these, we need some more notation.
The inverse Cartan matrix of the Lie algebra A$_{\ell-1}$ is denoted
by $B_{\ell-1}$, and $\vec{\mu}$ and $\vec{\eps}_j$ are
$(\ell-1)$-dimensional (column) vectors with entries
$\vec{\mu}_j=\mu_j$ and $(\vec{\epsilon}_j)_{m}=\delta_{j,m}$.
Furthermore, we need the $k$-dimensional vector
\begin{equation}
\vec{Q}_{i,i',\ell} = 
\vec{\e}_i + \vec{\e}_{i+2} + \cdots +
\vec{\e}_{i'} + \vec{\e}_{i'+2} + \cdots +
\vec{\e}_{\ell+1} + \vec{\e}_{\ell+3} + \cdots ,
\end{equation}
with $\vec{\e}_j = \vec{0}$ for $j \geq k+1$.
Using this notation, we are led to the following conjecture.
\begin{conjecture}
For all $k\geq 1$, $1\leq \ell \leq k$, $1\leq i \leq k+1$,
$1\leq i'\leq \ell+1$ and $|q|<1$, the $q$-series
{\rm
\begin{equation}
\renewcommand{\arraystretch}{0.7}
q^{(i'+i-2)/4} \hspace{-4mm}
\sum_{\begin{array}{c}
\sc m_1,m_2,\ldots,m_k\geq 0\\
\sc m_j \equiv (\vec{Q}_{i,i',\ell})_j \mod{2}
\end{array}}
\hspace{-4mm}
\frac{q^{\displaystyle \, \case{1}{4} \, \vec{m}^T C_k (\vec{m}
+ 2\,\vec{\e}_k - 2\,\vec{\e}_{i-1})}}{(q)_{m_{\ell}}}
\prod_{\begin{array}{c}
\sc j=1\\
\sc j \neq \ell
\end{array} }^k
\Bin{\case{1}{2}\bigl({\cal I}_k \, \vec{m}
+ \vec{\e}_{i-1}
+\vec{\e}_{i'-1}-\vec{\e}_k \bigr)}{m_j}
\label{Fdual}
\end{equation}}
equals
{\rm \begin{eqnarray}
\renewcommand{\arraystretch}{0.7}
\lefteqn{
q^{(k+r-s+1)(k-r-s+1)/(4\ell)}
\; \frac{1}{(q)_{\infty}} \;
\sum_{n=0}^{\ell-1}
\hspace{-2mm}
\sum_{\begin{array}{c}
\sc \mu_1,\ldots,\mu_{\ell-1}\geq 0 \\
\sc n + \ell (B_{\ell-1} \vec{\mu})_1
\equiv 0 \mod{\ell}
\end{array}}
\hspace{-4mm}
\frac{q^{\displaystyle \, \vec{\mu}^T B_{\ell-1} (\vec{\mu}
 - \vec{\epsilon}_r)}}{(q)_{\mu_1} \cdots (q)_{\mu_{\ell-1}}}} \nonumber \\
& &\renewcommand{\arraystretch}{0.7}
\mbox{} \times \biggl\{
\sum_{\begin{array}{c}
\sc j=-\infty \\
\sc n+(k-s-r+1)/2+(2k+3)j\equiv 0 \mod{\ell}
\end{array}}^{\infty}
q^{j\bigl((2k-2\ell+3)(2k+3)j+
(2k+3)(k-\ell+1)-(2k-2\ell+3)s\bigr)/\ell} \biggr. \nonumber \\
& &\renewcommand{\arraystretch}{0.7}
\biggl. \qquad
- \hspace{-4mm}
\sum_{\begin{array}{c}
\sc j=-\infty \\
\sc n+(k+s-r+1)/2+(2k+3)j\equiv 0 \mod{\ell}
\end{array}}^{\infty}
q^{\bigl((2k-2\ell+3)j+(k-\ell+1)\bigr)
\bigl((2k+3)j+s\bigr)/\ell}
\biggr\}
\label{rk}
\end{eqnarray}}
for $r\equiv k \mod{2}$, and equals
{\rm \begin{eqnarray}
\renewcommand{\arraystretch}{0.7}
\lefteqn{
q^{(k+r-s+2)(k-r-s+2)/(4\ell)}
\; \frac{1}{(q)_{\infty}} \; 
\sum_{n=0}^{\ell-1}
\sum_{\begin{array}{c}
\sc \mu_1,\ldots,\mu_{\ell-1}\geq 0 \\
\sc n + \ell (B_{\ell-1} \vec{\mu})_1
\equiv 0 \mod{\ell} \end{array}}
\frac{q^{\displaystyle \, \vec{\mu}^T B_{\ell-1} (\vec{\mu}
 - \vec{\epsilon}_r)}}{(q)_{\mu_1} \cdots (q)_{\mu_{\ell-1}}}} \nonumber \\
& &\renewcommand{\arraystretch}{0.7}
\mbox{} \times \biggl\{
\sum_{\begin{array}{c}
\sc j=-\infty \\
\sc n-(k-s+r+2)/2-(2k+3)j\equiv 0 \mod{\ell}
\end{array}}^{\infty}
q^{j\bigl((2k-2\ell+3)(2k+3)j+
(2k+3)(k-\ell+2)-(2k-2\ell+3)s\bigr)/\ell} \biggr. \nonumber \\
& &\renewcommand{\arraystretch}{0.7}
\biggl. \qquad
- \hspace{-4mm}
\sum_{\begin{array}{c}
\sc j=-\infty \\
\sc n-(k+s+r+2)/2-(2k+3)j\equiv 0 \mod{\ell}
\end{array}}^{\infty}
q^{\bigl((2k-2\ell+3)j+(k-\ell+2)\bigr)
\bigl((2k+3)j+s\bigr)/\ell}
\biggr\}
\label{rnk}
\end{eqnarray}}
for $r\not\equiv k \mod{2}$. 
\end{conjecture}
Since conjecture~\ref{conj} is proven for
$\ell=1$ and $k$, we can for these particular values
claim the above as theorem. 
In fact, for $\ell=1$, the above was first conjectured in ref.~\cite{KKMMb}
and proven in~\cite{FQ94}.
In refs.~\cite{KMQ,BNY} expressions for the branching functions of
the $({\rm A}^{(1)}_1)_{M} \times ({\rm A}^{(1)}_1)_{N} /
({\rm A}^{(1)}_1)_{M+N}$
coset conformal field theories were given similar 
to (\ref{rk}) and (\ref{rnk}).
This similarity suggests that (\ref{Fdual}), (\ref{rk}) and (\ref{rnk})
correspond to the branching functions of the coset 
$({\rm A}^{(1)}_1)_{\ell} \times ({\rm A}^{(1)}_1)_{k-\ell-1/2} /
({\rm A}^{(1)}_1)_{k-1/2}$ of fractional level.

A very last comment we wish to make is that
there exist other polynomial identities
than those discussed in this paper which
imply the Andrews--Gordon identity~\ref{An} and which
involve the $q$-multinomial coefficients.
\begin{theorem}\label{t7}
For all $k \leq 1$ and $1 \leq i \leq k+1$,
\begin{equation}
\sum_{a=0}^L \dMult{L}{a}{k-i+1}{k} =
\sum_{j=-L}^L
(-)^j
q^{j\bigl( (2k+3)(j+1)-2i\bigr)/2} \frac{(q)_L}{(q)_{L-j}(q)_{L+j}} \: .
\label{e7}
\end{equation}
\end{theorem}
Note that for $i=k+1$ the left-hand side is the generating
function of partitions with at most $k$ successive Durfee
squares and with largest part $\leq L$.

The proof of theorem~\ref{t7} follows readily using the
Bailey lattice of refs.~\cite{AAB}.
For $k=1$ (\ref{e7}) was first obtained by Rogers~\cite{Rogers94}.
For other $k$ it is implicit in refs.~\cite{AAB,Andrews84}.

\section*{Acknowledgements}
I thank Anne Schilling for helpful and stimulating
discussions on the $q$-multinomial coefficients.
Especially her communication of equation (\ref{frec})
has been indispensable for proving proposition~\ref{pqr2}.
I thank Alexander Berkovich for very constructive discussions 
on the nature of the fermi-gas of section~\ref{pt4}.
I wish to thank Professor G.~E.~Andrews for drawing my attention
to the relevance of equation (\ref{fk}) and
Barry McCoy for electronic lectures on
the history of the Rogers--Ramanujan identities.
Finally, helpful and interesting discussions with
Omar Foda and Peter Forrester are greatfully acknowledged.
This work is supported by the Australian Research Council.

\appendix

\nsection{Proof of $q$-multinomial relations}
In this section we prove the various claims 
concerning the $q$-multinomial coefficients
made in section~\ref{secqB}. 

Let us start proving the symmetry properties
(\ref{qs}) of lemma~\ref{lemsym}.
First we take the definition (\ref{qM})
and make the change variables $j_{\ell} \to
L-j_{k-\ell+1}$ for all $\ell=1,\ldots,k$.
This changes the restriction on the sum to 
$j_1+ \cdots +j_k=kL-a$, changes the exponent
of $q$ to
\begin{equation}
\sum_{\ell=1}^{k-1} (L-j_{k-\ell}) j_{k-\ell+1}
-\sum_{\ell=k-p}^{k-1} (L-j_{k-\ell}),
\label{qform}
\end{equation}
but leaves the product over the $q$-binomials invariant.
We now perform a simple rewriting of (\ref{qform}) as
follows
\begin{eqnarray}
(\mbox{\ref{qform}}) &=& \sum_{\ell=1}^{k-1}(L-j_{\ell})j_{\ell+1}
+ \sum_{\ell=1}^p j_{\ell} -pL \hspace{4cm}
(\mbox{by } \ell \to k-\ell) \nonumber \\
&=& 
\sum_{\ell=1}^{k-1}(L-j_{\ell})j_{\ell+1}
- \sum_{\ell=p}^{k-1} j_{\ell+1} +(k-p)L-a 
\hspace{1cm} (\mbox{by } j_1+\cdots + j_k=kL-a),
\end{eqnarray}
which proves the first claim of the lemma.
The second statement in the lemma follows 
for example, by noting that $\Mults{L}{a}{k}{k}=
q^{-a} \Mults{L}{a}{0}{k}$.

The proof of the tautologies (\ref{taut}) of proposition~\ref{ptaut}
is somewhat more involved and we proceed inductively.
For $L=0$ (\ref{taut}) is obviously correct, thanks to
\begin{equation}
\Mult{0}{a}{p}{k} = \delta_{a,0}.
\end{equation}
Now assume (\ref{taut}) holds true for all $L'=0,\ldots,L$.
To show that this implies (\ref{taut}) for $L'=L+1$,
we substitute the fundamental recurrence (\ref{frec})
into (\ref{taut}) with $L$ replaced by $L+1$.
After some cancellation of terms and division by $(1-q^{L+1})$,
this simplifies to
\begin{equation}
\sum_{m=0}^M q^{m L} \Mult{L}{a-m}{m}{k}=
\sum_{m=0}^M q^{m L} \Mult{L}{kL-a-m+M}{m}{k},
\label{t2}
\end{equation}
where we have replaced $k-p-1$ by $M$.
Since in (\ref{taut}) we have $p=-1,\ldots,k-1$,
(\ref{t2}) should hold for $M=0,\ldots,k$.
A set of equations equivalent to this is obtained
by taking $\mbox{(\ref{t2})}_{M=0}$ and 
$\mbox{(\ref{t2})}_M-\mbox{(\ref{t2})}_{M-1}$ for $M=1,\ldots,k$.
In formula this new set of equations reads
\begin{eqnarray}
\lefteqn{
\sum_{m=0}^{M-1} q^{m L}\left\{ \Mult{L}{kL-a-m+M-1}{m}{k} 
-q^L \Mult{L}{kL-a-m+M-1}{m+1}{k} \right\}}
\nonumber \\
& & \hspace{3cm} = 
\Mult{L}{kL-a+M}{0}{k}-q^{M L} \Mult{L}{a-M}{M}{k},
\end{eqnarray}
for $M=0,\ldots,k$.
Now we use the induction assumption on the term
within the curly braces, and the second symmetry
relation of (\ref{qs}) on the first term of the right-hand side.
This yields
\begin{equation}
\sum_{m=0}^{M-1} q^{m L} \left\{\Mult{L}{a-M}{m}{k} 
-q^L \Mult{L}{a-M}{m+1}{k} \right\}
=\Mult{L}{a-M}{0}{k}-q^{M L} \Mult{L}{a-M}{M}{k}.
\end{equation}
Expanding the sum, all but two terms on the left-hand side
cancel, yielding the right-hand side.

Finally we have to show equation~(\ref{qrec})
of proposition~\ref{pqr2} to be true.
We approach this problem indirectly and 
will in fact show that the right-hand side of (\ref{qrec})
can be transformed into the right-hand side
of (\ref{frec}) by multiple application of the
tautologies~(\ref{taut}) and the symmetries (\ref{qs}).
For the sake of convenience, we restrict our attention
to the case $k$ and $p$ even, and replace $L$ in (\ref{frec})
and (\ref{qrec}) by $L+1$. The other choices for the parity of
$k$ and $p$ follow in analogous manner, and the details will
be omitted.

Rewriting the right-hand side of (\ref{qrec}) by
replacing $L$ by $L+1$, using the even parity
of $k$ and $p$, and replacing $p$ by $k-M$, gives
\begin{eqnarray}
\lefteqn{
\renewcommand{\arraystretch}{0.6}
\sum_{\begin{array}{c}
\sc m=0 \\
\sc m \mbox{\scriptsize even}
\end{array} }^M
q^{mL}
\Mult{L}{a-\case{1}{2}(m+M)}{m}{k}
+
\sum_{\begin{array}{c}
\sc m=1 \\
\sc m \mbox{\scriptsize odd}
\end{array} }^{M-1}
q^{mL}
\Mult{L}{kL-a-\case{1}{2}(m-M+1)}{m}{k}}  \nonumber \\
& &
\renewcommand{\arraystretch}{0.6}
+\sum_{\begin{array}{c}
\sc m=M+2 \\
\sc m \mbox{\scriptsize even}
\end{array} }^k
q^{\case{1}{2}\bigl((2L+1)M-m\bigr)}
\Mult{L}{a-\case{1}{2}(m+M)}{m}{k}
\nonumber \\
& &
\renewcommand{\arraystretch}{0.6}
+ \sum_{\begin{array}{c}
\sc m=M+1 \\
\sc m \mbox{\scriptsize odd}
\end{array} }^{k-1}
q^{k(L+1)+\case{1}{2}\bigl((2L+3)M-m+1\bigr)-2a}
\Mult{L}{k(L+1)-a-\case{1}{2}(m-M-1)}{m}{k}.
\label{eq1}
\end{eqnarray}

The proof that this equals the right-hand side of (\ref{frec})
(with $L$ replaced by $L+1$ and $p$ by $k-M$) 
breaks up into two independent
steps, both of which will be given as a lemma.
First, we have
\begin{lemma}\label{p1}
The top-line of equation (\ref{eq1}) equals
\begin{equation}
\sum_{m=0}^M q^{mL} \Mult{L}{a-m}{m}{k}.
\label{ep1}
\end{equation}
\end{lemma}
Second,
\begin{lemma}\label{p2}
The bottom-two lines of equation (\ref{eq1}) equal
\begin{equation}
\sum_{m=M+1}^k
q^{(L+1)M-m}
\Mult{L}{a-m}{m}{k}.
\label{ep2}
\end{equation}
\end{lemma}
Clearly, application of these two lemmas
immediately establishes the wanted result.

At the core of the proof of both lemmas is yet
another result, which can be stated as
\begin{lemma}\label{p3}
For $M$ even and $\ell=0,\ldots,\case{1}{2}M$, 
the following function is independent of $\ell$:
\begin{eqnarray}
\lefteqn{
F_{\ell}(M,a)=
\sum_{m=M-\ell}^M q^{mL} \Mult{L}{a-m}{m}{k}
+ \sum_{m=0}^{\ell-1} q^{mL} \Mult{L}{kL-a-m+\case{1}{2}M-1}{m}{k}} \\
&&+ \!\!
\renewcommand{\arraystretch}{0.6}
\sum_{\begin{array}{c}
\sc m=\ell \\
\sc m \equiv \ell \mod{2}
\end{array} }^{M-\ell-2} \!\!
q^{mL} \left\{
\Mult{L}{a-\case{1}{2}(m+M-\ell)}{m}{k}
+ q^L
\Mult{L}{kL-a-\case{1}{2}(m-M+\ell)-1}{m+1}{k} \right\}.
\nonumber
\end{eqnarray}
\end{lemma}
The proof of this is simple. First we apply the tautology
(\ref{taut}) to the term within the curly braces, yielding
\begin{eqnarray}
\lefteqn{
F_{\ell}(M,a)=
\sum_{m=M-\ell}^M q^{mL} \Mult{L}{a-m}{m}{k}
+ \sum_{m=0}^{\ell-1} q^{mL} \Mult{L}{kL-a-m+\case{1}{2}M-1}{m}{k}} \\
&&+ \!\!
\renewcommand{\arraystretch}{0.6}
\sum_{\begin{array}{c}
\sc m=\ell \\
\sc m \equiv \ell \mod{2}
\end{array} }^{M-\ell-2} \!\!
q^{mL} \left\{
\Mult{L}{kL-a-\case{1}{2}(m-M+\ell)-1}{m}{k} 
+ q^L
\Mult{L}{a-\case{1}{2}(m+M-\ell)}{m+1}{k} \right\}.
\nonumber
\end{eqnarray}
After separating the $m=\ell$ term in the first 
and the $m=M-\ell-2$ term in the second term within 
the curly braces, 
we change the summation variable $m\to m-2$ in the sum
over the second term within the braces.
This results in
\begin{eqnarray}
\lefteqn{
F_{\ell}(M,a)=
\sum_{m=M-\ell-1}^M q^{mL} \Mult{L}{a-m}{m}{k}
+ \sum_{m=0}^{\ell} q^{mL} \Mult{L}{kL-a-m+\case{1}{2}M-1}{m}{k}} 
\nonumber \\
&&\quad + \!\!
\renewcommand{\arraystretch}{0.6}
\sum_{\begin{array}{c}
\sc m=\ell+1 \\
\sc m \equiv \ell+1 \mod{2}
\end{array} }^{M-\ell-3} \!\!
q^{mL} \left\{
\Mult{L}{a-\case{1}{2}(m+M-\ell-1)}{m}{k} \right.
\nonumber \\[-3mm]
&&\hspace{3.5cm} +  \left. q^L
\Mult{L}{kL-a-\case{1}{2}(m-M+\ell+1)-1}{m+1}{k} \right\} = 
F_{\ell+1}(M,a).
\end{eqnarray}

The proof of the lemmas~\ref{p1} and \ref{p2} readily follows
from lemma~\ref{p3}.
To prove lemma~\ref{p1}, note that the 
top-line of (\ref{eq1}) is nothing but $F_0(M)$.
Since this is equal to $F_{\frac{1}{2}M}(M)$,
we get
\begin{equation}
\mbox{top-line of (\ref{eq1})}=
\sum_{m=\frac{1}{2}M}^M q^{mL} \Mult{L}{a-m}{m}{k}
+ \sum_{m=0}^{\frac{1}{2}M-1} q^{mL} 
\Mult{L}{kL-a-m+\case{1}{2}M-1}{m}{k}.
\end{equation}
Applying equation (\ref{t2}) with $M$ replaced by $\frac{1}{2}M-1$,
to the the second sum, we simplify to equation (\ref{ep1}) thus 
proving our lemma.

To prove lemma~\ref{p2}, we apply the first
symmetry relation of (\ref{qs}) to
all $q$-multinomials in the bottom-two
lines of (\ref{eq1}).  After changing $m\to k-m$ in the
second line and $m\to k-m-1$ in the third line, the
last two lines of (\ref{eq1}) combine to
\begin{equation}
\renewcommand{\arraystretch}{0.6}
f_a
\sum_{\begin{array}{c}
\sc m=0 \\
\sc m \mbox{\scriptsize even}
\end{array} }^{k-M-2}
q^{mL} \left\{
\Mult{L}{kL-a-\case{1}{2}(m-M-k)}{m}{k}
+
q^L \Mult{L}{a-\case{1}{2}(m+M+k)-1}{m+1}{k}
\right\},
\label{sim}
\end{equation}
with $f_a =q^{(L+1)M-a}$.
This we recognize as 
\begin{equation}
f_a \left\{
F_0(k-M,kL+k-a)-q^{(k-M)L}\Mult{L}{kL-a+M}{k-M}{k}\right\}.
\end{equation}
Replacing the first term by $F_{\frac{1}{2}(k-M)}(k-M,kL+k-a)$, gives
\begin{equation}
f_a
\sum_{m=0}^{\frac{1}{2}(k-M)-1}
q^{mL}
\Mult{L}{a-\case{1}{2}(M+k)-m-1}{m}{k}
+f_a \sum_{m=\frac{1}{2}(k-M)}^{k-M-1}
q^{mL}
\Mult{L}{kL-a-m+k}{m}{k} .
\end{equation}
Applying equation (\ref{t2}) with $M$ replaced by $\frac{1}{2}(k-M)-1$
and $a$ by $a-\frac{1}{2}(M+k)-1$, to the first sum, this simplifies to
\begin{equation}
f_a \sum_{m=0}^{k-M-1}
q^{mL}
\Mult{L}{kL-a-m+k}{m}{k} .
\end{equation}
Finally using the symmetry (\ref{qs}), recalling the
definition of $f_a$ and changing $m\to k-m$, we get
equation (\ref{ep2}).

\end{document}